\documentclass[a4paper, 12pt, leqno]{amsart}
\usepackage{amsmath,amsfonts,amssymb}
\usepackage{mathrsfs}
\usepackage{amscd}
\usepackage{amsthm}
\usepackage{mathrsfs}
\usepackage{hyperref}
\usepackage{epic}
\usepackage{mathtools}

\title{The boundness of weighted Coxeter groups of rank 3}
\author{Jianwei Gao}
\keywords{Weighted Coxeter group, Hecke algebra, Two-sided cell, Left cell, A-function}

\address{Jianwei Gao\\
Beijing International Center for Mathematical Research\\
Peking University\\
No.5 Yiheyuan Road\\
Beijing, 100871\\
People's Republic of China}
\email{gaojianwei@bicmr.pku.edu.cn}

\begin{document}

\maketitle

\begin{abstract}
We prove that a weighted Coxeter group $(W,S,L)$ is bounded in the sense of G.Lusztig if the rank of $W$ is 3.
\end{abstract}

\setcounter{section}{-1}
\section{Introduction}

\medskip

Let $(W,S,L)$ be a weighted Coxeter group. In [L3, 13.4], G.Lusztig conjectured that the maximal weight value of the longest elements of the finite parabolic subgroups of $W$ is a bound for $(W,S,L)$. This property is referred as boundness of a weighted Coxeter group([L3, 13.2]). When $W$ is finite, this conjecture is clear. In [L1, 7.2], G.Lusztig proved this conjecture when $W$ is an affine Weyl group and $L=l$, but the same proof remains valid without the assumption $L=l$, see [L3, Chapter 13]. In [SY, 3.2], J.Shi and G.Yang proved this conjecture when $W$ has complete Coxeter graph. In [Zhou, 2.1], P.Zhou proved this conjecture when $rank(W)=3$ and $L=l$. In this paper, we will prove this conjecture when $rank(W)=3$ without the assumption $L=l$, see Theorem 2.1. Then we can get some interesting consequences and describe the lowest two-sided cell of $W$ in this case. The author would like to thank N.Xi for his comments.

\medskip

\section{Preliminaries}

\medskip

\noindent{\bf 1.1.}
In this paper, for any Coxeter group $(W,S)$, we assume the generating set $S$ is finite. We call $|S|$ the rank of $(W,S)$ and denote it by $rank(W)$. We use $l$ for the length function and $\leq$ for the Bruhat order of $W$. The neutral element of $W$ will be denoted by $e$. For $x\in W$, we set $$\mathcal{L}(x)=\{s\in S|sx<x\},\ \mathcal{R}(x)=\{s\in S|xs<x\}.$$
For $s,t \in W$, let $m_{st}\in\mathbb{Z}_{\geq1}\bigcup\{\infty\}$ be the order of $st$ in $W$.

For any $I\subseteq S$, let $W_{I}=\langle I\rangle$. Then $(W_{I},I)$ is also a Coxeter group, called a parabolic subgroup of $(W,S)$. Denote the longest element of $W_{I}$ by $w_{I}$ if $|W_{I}|<\infty$. For $s,t\in S$, $s\neq t$, we use $W_{st}$ instead of $W_{\{s,t\}}$ and $w_{st}$ instead of $w_{\{s,t\}}$.

For $w_{1},w_{2},\cdots,w_{n}\in W$, we often use the notation $w_{1}\cdot w_{2}\cdot\  \cdots \ \cdot w_{n}$ instead of $w_{1}w_{2}\cdots w_{n}$ if
$l(w_{1}w_{2}\cdots w_{n})=l(w_{1})+l(w_{2})+\cdots +l(w_{n})$.

\medskip

\noindent{\bf 1.2.}
Let $(W,S)$ be a Coxeter group. A map $L:W\longrightarrow \mathbb{Z}$ is called a weight function if $L(ww')=L(w)+L(w')$ for any $w,w'\in W$ with $l(ww')=l(w)+l(w')$. Then we call $(W,S,L)$ a weighted Coxeter group. For any $J\subseteq S$,
it is obvious that the parabolic subgroup $(W_{J},J,L|_{W_{J}})$ is also a weighted Coxeter group.

In this paper, the weight function $L$ for any weighted Coxeter group $(W,S,L)$ is assumed to be positive, that is, $L(s)>0$ for any $s\in S$.

\medskip

\noindent{\bf 1.3.}
Let $(W,S,L)$ be a weighted Coxeter group and $\mathbb{Z}[v,v^{-1}]$ be the ring of Laurent polynomials in an indeterminate $v$ with integer coefficients. For $f=\sum \limits_{n\in \mathbb{Z}}a_{n}v^{n}\in \mathbb{Z}[v,v^{-1}]\setminus\{0\}$, we define $deg\ f=\max \limits_{n\in \mathbb{Z}\atop a_{n}\neq 0} n$. Complementally, we define $deg\ 0=-\infty$.

For $w\in W$, set $v_{w}=v^{L(w)}\in \mathbb{Z}[v,v^{-1}]$.
The Hecke algebra $\mathcal{H}$ of $(W,S,L)$ is the unital associative $\mathbb{Z}[v,v^{-1}]-$algebra defined by the generators $T_{s}(s\in S)$ and the relations
$$
(T_{s}-v_{s})(T_{s}+v_{s}^{-1})=0,\ \forall s\in S.
$$
$$
\underbrace{T_{s}T_{t}T_{s}\cdots}_{m_{st}\ factors}\ =\ \underbrace{T_{t}T_{s}T_{t}\cdots}_{m_{st}\ factors},\ \forall s,t\in S,\ m_{st}<\infty.
$$

Obviously, $T_{e}$ is the multiplicative unit of $\mathcal{H}$. For any $w\in W$, we define $T_{w}=T_{s_{1}}T_{s_{2}}\cdots T_{s_{n}}\in \mathcal{H}$, where $w=s_{1}s_{2}\cdots s_{n}$ is a reduced expression of $w$ in $W$. Then $T_{w}$
is independent of the choice of reduced expression and $\{T_{w}|w\in W\}$ is a $\mathbb{Z}[v,v^{-1}]-$basis of $\mathcal{H}$, called the standard basis.
We define $f_{x,y,z}\in \mathbb{Z}[v,v^{-1}]$ for any $x,y,z\in W$ by the identity
$$
T_{x}T_{y}=\sum \limits_{z\in W} f_{x,y,z}T_{z}.
$$
The following involutive automorphism of rings is useful, called the bar involution:
$$\bar{}:\mathcal{H}\longrightarrow \mathcal{H}\ \ \ \ \ $$
$$v^{n}\longrightarrow v^{-n}.$$
$$T_{s}\longrightarrow T_{s}^{-1}.$$
We have $\overline{T_{w}}=T_{w^{-1}}^{-1}$ for any $w\in W$. We set
$$\mathcal{H}_{\leq 0}=\bigoplus \limits_{w\in W}\mathbb{Z}[v^{-1}]T_{w},\ \mathcal{H}_{< 0}=\bigoplus \limits_{w\in W}v^{-1}\mathbb{Z}[v^{-1}]T_{w}.$$
We can get the following facts by easy computation.

\medskip

\noindent{\bf Lemma 1.4.}
(1) For any $x,y,\in W$, we have
$$f_{x,y,e}=\delta_{x,y^{-1}}.$$
(2) For any $x,y,z\in W$, we have
$$deg\ f_{x,y,z}\leq min\{L(x),L(y),L(z)\}.$$
(3) For any finite parabolic subgroup $W_{I}$ of $W$, $x\in W_{I}$, we have
$$deg\ f_{w_{I},w_{I},x}=L(x).$$
(4) Assume $s,t\in S$, $m_{st}<\infty$, $x,y,z\in W_{st}$ such that $deg\ f_{x,y,z}=L(z)$. If $l(z)\geq 2$,
then $x=y=w_{st}$. If $z=s$, then $y=x^{-1}$, $s\in \mathcal{L}(x)$, $s\in \mathcal{R}(y)$.

\medskip

Define the degree map
$$deg:\ \mathcal{H}\longrightarrow \mathbb{Z}\cup\{-\infty\}\ \ \ \ \ \ \ \ \ $$
$$\ \ \ \ \sum_{w\in W}f_{w}T_{w}\longrightarrow Max\{deg\ f_{w}|w\in W\}.$$
And we set $N=\max \limits_{I\subseteq S\atop |W_{I}|<\infty} L(w_{I})$. G.Lusztig gave the following conjecture in [L3, 13.4].

\medskip

\noindent{\bf Conjecture 1.5.}
Let $(W,S,L)$ be a weighted Coxeter group, $S$ is finite and $L$ is positive, then $N$ is a bound for $(W,S,L)$. Namely, $deg\ (T_{x}T_{y})\leq N$ for all $x,y \in W$.

\medskip

\noindent{\bf Remark 1.6.}
This conjecture is very important in studying the cells of weighted Coxeter groups. For example, if this conjecture is true,
then the a-function is also bounded by this $N$. Further, $W$ has a lowest two-sided cell, see section 7. We can also consider conjectures P1-P15 and the ring $J$, see [L3, Chapter 14, Chapter 18].

When $W$ is a finite Coxeter group, this conjecture can be proved using Lemma 1.4(2). In [L1, 7.2], G.Lusztig proved this conjecture when $W$ is an affine Weyl group. In [SY, 3.2], J.Shi and G.Yang proved this conjecture when $W$ has complete Coxeter graph. In [Zhou, 2.1], P.Zhou proved this conjecture when $rank(W)=3$ and $L=l$, so this conjecture is also true when $rank(W)=3$ and $L=nl$ for some $n\in \mathbb{Z}_{\geq1}$. In this paper, we will prove this conjecture when $rank(W)=3$ and $L$ is positive, see Theorem 2.1.

\medskip

\noindent{\bf 1.7.}
For any $w\in W$, there exists a unique element $C_{w}\in \mathcal{H}_{\leq 0}$ such that
$\overline{C_{w}}=C_{w}$ and $C_{w}-T_{w}\in \mathcal{H}_{< 0}$. The elements $\{C_{w}|w\in W\}$ form a $\mathbb{Z}[v,v^{-1}]-$basis of $\mathcal{H}$, called the Kazhdan-Lusztig basis. We define $h_{x,y,z}\in \mathbb{Z}[v,v^{-1}]$ for any $x,y,z\in W$ such that
$$
C_{x}C_{y}=\sum \limits_{z\in W} h_{x,y,z}C_{z}.
$$

Using Kazhdan-Lusztig basis, we can define the preorders $\underset{L}{\leq}$, $\underset{R}{\leq}$, $\underset{LR}{\leq}$ on $W$. These preorders give rise to equivalence relations $\underset{L}{\sim}$, $\underset{R}{\sim}$, $\underset{LR}{\sim}$ on $W$ respectively. The equicalence classes are called left cells, right cells and two-sided cells of $W$. Then we have partial orders $\underset{L}{\leq}$, $\underset{R}{\leq}$, $\underset{LR}{\leq}$ on the sets of left cells, right cells and two-sided cells of $W$ respectively. For $x,y\in W$, we have $x\underset{L}{\leq}y$ if and only if $x^{-1}\underset{R}{\leq}y^{-1}$.

\medskip

Now we assume Conjecture 1.5 holds, so we have $deg\ h_{x,y,z}\leq N$ for any $x,y,z\in W$. Then we can define the a-function
$$a:W\longrightarrow \mathbb{N}\ \ \ \ \ \ \ \ \ \ \ $$
$$\ \ \ \ \ \ \ \ \ \ \ w\longrightarrow \max \limits_{x,y\in W} deg\ h_{x,y,w}.$$

\medskip

For any $x,y,z\in W$, we define $\beta_{x,y,z},\gamma_{x,y,z}\in \mathbb{Z}$ such that
$$f_{x,y,z^{-1}}=\beta_{x,y,z}v^{N}+\mbox{lower degree terms}.$$
$$h_{x,y,z^{-1}}=\gamma_{x,y,z}v^{a(z)}+\mbox{lower degree terms}.$$
\noindent{\bf Lemma 1.8.} Let $x,y,z\in W$.
\\(1) We have $\beta_{x,y,z}=\beta_{y,z,x}=\beta_{z,x,y}$.
\\(2) We have $\gamma_{x,y,z}=\gamma_{y^{-1},x^{-1},z^{-1}}$.
\\(3) If $\beta_{x,y,z}\neq 0$, then $x\underset{L}{\sim}y^{-1}$, $y\underset{L}{\sim}z^{-1}$, $z\underset{L}{\sim}x^{-1}$, $a(x)=a(y)=a(z)=N$, and $\beta_{x,y,z}=\gamma_{x,y,z}=\gamma_{y,z,x}=\gamma_{z,x,y}$.
\\(4) If $\gamma_{x,y,z}\neq 0$ and $a(z)=N$, then $\beta_{x,y,z}=\gamma_{x,y,z}\neq 0$.

\medskip

The detail of 1.7 and the proof of Lemma 1.8 can be found in [L3].

\medskip

\section{Weighted Coxeter Groups of Rank 3}

\medskip

From now on, we assume $(W,S,L)$ is a weighted Coxeter group of rank 3 and $L$ is positive. We set $N=\max \limits_{I\subseteq S\atop |W_{I}|<\infty} L(w_{I})$. The main result of this paper is the following.

\medskip

\noindent{\bf Theorem 2.1.}
We have $deg\ (T_{x}T_{y})\leq N$ for all $x,y \in W$. In other word, Conjecture 1.5 holds in this case.

\medskip

Let $S=\{r,s,t\}$. Conjecture 1.5 has been proved if $W$ is a finite Coxeter group or affine Weyl group or $W$ has complete Coxeter graph, so we may assume that $W$ is infinite, $m_{rt}=2$, $m_{sr}\geq m_{st}$.

\medskip

When $m_{sr}=\infty$, $m_{st}=2$, the Coxeter graph of $(W,S)$ is not connected. In this case, $(W,S)$ is the direct product of $(W_{sr}, \{s,r\})$ and $(\langle t\rangle, \{t\})$ as a Coxeter group. The parabolic subgroup $W_{sr}$ is an affine Weyl group of type $\widetilde{A_{1}}$ and $\langle t\rangle$ is a finite Coxeter group of two elements, Conjecture 1.5 holds in these two cases.

\medskip

When $m_{sr}=m_{st}=\infty$, the case is also simple. For any $x\in W$, we have $\mathcal{L}(x)=\{s\}$ or $\mathcal{L}(x)\subseteq \{r,t\}$, we also have $\mathcal{R}(x)=\{s\}$ or $\mathcal{R}(x)\subseteq \{r,t\}$. For any $x,y\in W$, if $\mathcal{R}(x)=\{s\}$ and $\mathcal{L}(y)\subseteq \{r,t\}$, or $\mathcal{R}(x)\subseteq \{r,t\}$ and $\mathcal{L}(y)=\{s\}$, then we have $T_{x}T_{y}=T_{xy}$. Therefore, we can get $deg\ (T_{x}T_{y})\leq Max\{L(rt),L(s)\}$ for any $x,y\in W$ by easy computation.

\medskip

Summarizing the discussions above, we only need to consider the following three cases. We will deal with them in sections 3 to 5.
\\Case 1:  $m_{sr}=\infty>m_{st}\geq 3$.
\\
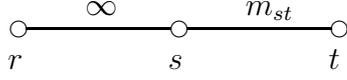
\begin{picture}(60,10)(-45,10)
\put(0,0){\circle{6}}
\put(60,0){\circle{6}}
\put(120,0){\circle{6}}
\put(-4,-15){$r$}
\put(56,-15){$s$}
\put(116,-15){$t$}
\put(3,0){\line(1,0){54}}
\put(63,0){\line(1,0){54}}
\put(25,5){$\infty$}
\put(85,5){$m_{st}$}
\end{picture}
\\\\\\Case 2:  $\infty>m_{sr}\geq m_{st}\geq 4$, $m_{sr}\geq5$.
\\
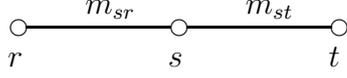
\begin{picture}(60,10)(-45,10)
\put(0,0){\circle{6}}
\put(60,0){\circle{6}}
\put(120,0){\circle{6}}
\put(-4,-15){$r$}
\put(56,-15){$s$}
\put(116,-15){$t$}
\put(3,0){\line(1,0){54}}
\put(63,0){\line(1,0){54}}
\put(25,5){$m_{sr}$}
\put(85,5){$m_{st}$}
\end{picture}
\\\\\\Case 3:  $\infty>m_{sr}\geq7$, $m_{st}=3$.
\\
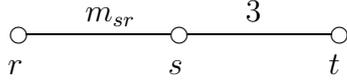
\begin{picture}(60,10)(-45,10)
\put(0,0){\circle{6}}
\put(60,0){\circle{6}}
\put(120,0){\circle{6}}
\put(-4,-15){$r$}
\put(56,-15){$s$}
\put(116,-15){$t$}
\put(3,0){\line(1,0){54}}
\put(63,0){\line(1,0){54}}
\put(25,5){$m_{sr}$}
\put(85,5){$3$}
\end{picture}

\vspace{11 mm}

\section{Case 1}

\medskip

In this section, we prove Theorem 2.1 for the case $m_{sr}=\infty>m_{st}\geq 3$. In this case, $N=Max\{L(rt),L(w_{st})\}$. It is easy to see

\medskip

\noindent{\bf Lemma 3.1.}
Let $x\in W$.
\\(1) If $s\in \mathcal{R}(x)$, then $r\notin \mathcal{R}(x)$.
\\(2) If $s\in \mathcal{L}(x)$, then $r\notin \mathcal{L}(x)$.
\\(3) If $x=x_{1}\cdot st$, then $r\notin \mathcal{R}(x)$.
\\(4) If $x=ts\cdot x_{1}$, then $r\notin \mathcal{L}(x)$.
\\(5) If $x=x_{1}\cdot rs$, then $\mathcal{R}(x)=\{s\}$.
\\(6) If $x=sr\cdot x_{1}$, then $\mathcal{L}(x)=\{s\}$.

\medskip

\noindent{\bf Lemma 3.2.}
Let $x,y\in W$.
\\(1) If $w\in W_{sr}$, $l(w)\geq 4$, $\mathcal{R}(x),\ \mathcal{L}(y)\subseteq\{t\}$, then $l(xwy)=l(x)+l(w)+l(y)$.
\\(2) If $\mathcal{R}(x),\ \mathcal{L}(y)\subseteq\{s\}$, then $l(xtry)=l(xrty)=l(x)+l(y)+2$.
\\(3) If $w\in W_{st}$, $l(w)\geq 2$, $\mathcal{R}(x),\ \mathcal{L}(y)\subseteq\{r\}$, then $l(xwy)=l(x)+l(w)+l(y)$.

\medskip

\noindent{\bf Proof.}
(1) See [Zhou, 5.7].
\\(2) See [Zhou, 5.9].
\\(3) When $l(w)\geq 3$, see [Zhou, 5.8]. We only prove the case of $w=st$ because the case of $w=ts$ is similar.

If $l(x)\leq 1$ or $l(y)\leq 1$, by Lemma 3.1(2),(3), we have $l(xwy)=l(x)+l(w)+l(y)$.

If $l(x)\geq2$ and $l(y)\geq2$, we may assume $x=x_{1}\cdot sr$, $y=rs\cdot y_{1}$ for some $x_{1},y_{1}\in W$. Since $\mathcal{L}(y)\subseteq \{r\}$, we have $\mathcal{L}(sy_{1})=\{s\}$. By Lemma 3.1(5), we have $\mathcal{R}(x_{1}srs)=\{s\}$. So
$$
\begin{aligned}
l(xwy)&=l(x_{1}srstrsy_{1})\\
&=l(x_{1}srs)+l(tr)+l(sy_{1})\\
&=l(x_{1}sr)+l(st)+l(rsy_{1})\\
&=l(x)+l(w)+l(y)
\end{aligned}
$$
by Lemma 3.2(2).\hfill $\square$

\medskip

\noindent{\bf Proposition 3.3.}
We have $deg\ (T_{x}T_{y})\leq N$ for all $x,y \in W$ in case 1.

\noindent{\bf Proof.}
Now we prove Theorem 2.1 for this case. We use induction on $l(y)$. When $l(y)=0,1$, the conclusion is clear by Lemma 1.4(2). Now assume $l(y)\geq 2$
and the conclusion is true for $y'$ if $l(y')<l(y)$. We assume $x=x_{1}\cdot w$, $y=u\cdot y_{1}$, $w,u\in W_{st}$, $\mathcal{R}(x_{1}),\mathcal{L}(y_{1})\subseteq \{r\}$. Then we have
$$
\begin{aligned}
T_{x}T_{y}&=T_{x_{1}}(T_{w}T_{u})T_{y_{1}}\\
&=T_{x_{1}}(\sum_{p\in W_{st}}f_{w,u,p}T_{p})T_{y_{1}}\\
&=\sum_{p\in W_{st}}f_{w,u,p}T_{x_{1}p}T_{y_{1}}.
\end{aligned}
$$

We will prove $deg\ (f_{w,u,p}T_{x_{1}p}T_{y_{1}})\leq N$ for all $p\in W_{st}$. We have 7 cases.
\\$(1)\ p=e,\ l(x_{1})\leq 1$ or $l(y_{1})\leq 1$.

We have $deg\ (T_{x_{1}}T_{y_{1}})\leq L(r)$ and $f_{w,u,e}\in \mathbb{Z}$, so $deg\ (f_{w,u,e}T_{x_{1}}T_{y_{1}})\leq L(r)$.
\\$(2)\ p=e,\ l(x_{1})\geq2,\ l(x_{2})\geq2$.

Assume $x_{1}=x_{2}\cdot sr$, $y_{1}=rs\cdot y_{2}$ for some $x_{2},y_{2}\in W$. Then
$$
\begin{aligned}
T_{x_{1}}T_{y_{1}}&=T_{x_{2}sr}T_{rsy_{2}}\\
&=(v_{r}-v_{r}^{-1})T_{x_{2}srs}T_{y_{2}}+T_{x_{2}s}T_{sy_{2}}\\
&=(v_{r}-v_{r}^{-1})T_{x_{2}srsy_{2}}+T_{x_{2}s}T_{sy_{2}}.
\end{aligned}
$$
Since $l(sy_{2})<l(y)$, by induction hypothesis, we have $deg\ (T_{x_{2}s}T_{sy_{2}})\leq N$, so $deg\ (T_{x_{1}}T_{y_{1}})\leq N$.
Since $f_{w,u,e}\in \mathbb{Z}$, we get $deg\ (f_{w,u,e}T_{x_{1}}T_{y_{1}})\leq N$.
\\$(3)\ p=s,\ l(x_{1})\leq 1$ or $l(y_{1})\leq 1$.

By Lemma 3.1(1)(2), $T_{x_{1}s}T_{y_{1}}=T_{x_{1}sy_{1}}$. By Lemma 1.4(2), we have $deg\ f_{w,u,s}\leq L(s)$, so $deg\ (f_{w,u,s}T_{x_{1}s}T_{y_{1}})\leq L(s)$.
\\$(4)\ p=s,\ l(x_{1})\geq2,\ l(x_{2})\geq2$.

Assume $x_{1}=x_{2}\cdot sr$, $y_{1}=rs\cdot y_{2}$ for some $x_{2},y_{2}\in W$. By Lemma 3.2(1), $T_{x_{1}s}T_{y_{1}}=T_{x_{2}srs}T_{rsy_{2}}=T_{x_{2}srsrsy_{2}}$.
Since $deg\ f_{w,u,s}\leq L(s)$, we get $deg\ (f_{w,u,s}T_{x_{1}s}T_{y_{1}})\leq L(s)$.
\\$(5)\ p=t,\ l(x_{1})\leq 1$ or $l(y_{1})\leq 1$.

We have $deg\ (T_{x_{1}t}T_{y_{1}})\leq L(r)$ and $deg\ f_{w,u,t}\leq L(t)$, so we get $deg\ (f_{w,u,t}T_{x_{1}t}T_{y_{1}})$ $\leq L(rt)$.
\\$(6)\ p=t,\ l(x_{1})\geq2,\ l(x_{2})\geq2$.

Assume $x_{1}=x_{2}\cdot sr$, $y_{1}=rs\cdot y_{2}$ for some $x_{2},y_{2}\in W$. Since $\mathcal{R}(x_{1})\subseteq \{r\}$, $\mathcal{L}(y_{1})\subseteq \{r\}$, we have $\mathcal{R}(x_{2}s)=\mathcal{L}(sy_{2})=\{s\}$. Then by Lemma 3.2(2),
$$
\begin{aligned}
T_{x_{1}t}T_{y_{1}}&=(v_{r}-v_{r}^{-1})T_{x_{2}srt}T_{sy_{2}}+T_{x_{2}st}T_{sy_{2}}\\
&=(v_{r}-v_{r}^{-1})T_{x_{2}srtsy_{2}}+T_{x_{2}st}T_{sy_{2}}.
\end{aligned}
$$
\textcircled{\small{1}} $u=e$.

Since $l(sy_{2})<l(y)$, by induction hypothesis, we have $deg\ (T_{x_{2}st}T_{sy_{2}})\leq N$, thus $deg\ (T_{x_{1}t}T_{y_{1}})\leq N$ and $deg\ (f_{w,u,t}T_{x_{1}t}T_{y_{1}})\leq N$.
\\\textcircled{\small{2}} $u\neq e$.

Since $T_{x_{2}st}T_{tsy_{2}}=T_{x_{2}s}T_{sy_{2}}+(v_{t}-v_{t}^{-1})T_{x_{2}st}T_{sy_{2}}$, $l(tsy_{2})<l(y)$, $l(sy_{2})<l(y)$,
by induction hypothesis, we have $deg\ (T_{x_{2}st}T_{sy_{2}})\leq N-L(t)$, thus $deg\ (T_{x_{1}t}T_{y_{1}})\leq N-L(t)$. Since $deg\ f_{w,u,t}\leq L(t)$, we get $deg\ (f_{w,u,t}T_{x_{1}t}T_{y_{1}})\leq N$.
\\$(7)\ l(p)\geq 2$.

By Lemma 3.2(3), $T_{x_{1}p}T_{y_{1}}=T_{x_{1}py_{1}}$. Since $deg\ f_{w,u,p}\leq L(p)$, we get $deg\ (f_{w,u,p}T_{x_{1}p}T_{y_{1}})\leq L(w_{st})\leq N$.\hfill $\square$

\medskip

\section{Case 2}

\medskip

In this section, we prove Theorem 2.1 for the case $\infty>m_{sr}\geq m_{st}\geq 4$, $m_{sr}\geq5$. In this case, $N=Max\{L(w_{sr}),L(w_{st})\}>L(rt)$. Note that $L(srst)<N$. First we have

\medskip

\noindent{\bf Lemma 4.1.}
Let $w\in W$.
\\(1) If $w=w_{1}\cdot ts$, then $r\notin \mathcal{R}(w)$.
\\(2) If $w=w_{1}\cdot rs$, then $t\notin \mathcal{R}(w)$.
\\(3) If $w=w_{1}\cdot st,\ \mathcal{R}(w_{1}s)=\{s\}$, then $r\notin \mathcal{R}(w)$.
\\(4) If $w=w_{1}\cdot sr,\ \mathcal{R}(w_{1}s)=\{s\}$, then $t\notin \mathcal{R}(w)$.
\\(5) If $w=w_{1}\cdot tst$, then $r\notin \mathcal{R}(w)$.
\\(6) If $w=w_{1}\cdot rsr$, then $t\notin \mathcal{R}(w)$.
\\(7) There is no $w_{1},w_{2}\in W$ such that $w=w_{1}\cdot st=w_{2}\cdot sr$.
\\(8) If $\mathcal{L}(w)\subseteq\{r\}$, then $\mathcal{L}(rtw_{st}w)=\{r\}$.
\\(9) If $\mathcal{L}(w)\subseteq\{t\}$, then $\mathcal{L}(trw_{sr}w)=\{t\}$.

\medskip

\noindent{\bf Proof.}
We use induction on $l(w)$ to prove (1) and (2) simultaneously. When $l(w)=0,1,2$, (1) and (2) are clear. Now assume $l(w)\geq 3$ and (1) and (2) are true for $w'$
if $l(w')<l(w)$. If $w=w_{2}\cdot r=w_{1}\cdot ts$, we have $w=w_{3}\cdot w_{sr}$, so $w_{2}=w_{3}\cdot w_{sr}r,\ w_{1}t=w_{3}\cdot w_{sr}s$. We get $w_{1}t=w_{4}\cdot tr$.
So $w_{4}\cdot t=w_{1}tr=w_{3}\cdot w_{sr}sr$, contradict the induction hypothesis. If $w=w_{2}\cdot t=w_{1}\cdot rs$, we can find a contradiction similarly.
\\(3) If $r\in \mathcal{R}(w)$, then $\mathcal{R}(w)=\{r,t\},\ \mathcal{R}(w_{1}s)=\{s\}$, so any reduced expression of $w$
is ended by $st$, a contradiction.
\\(4) Similar to (3).
\\(5) By (1), $\mathcal{R}(w_{1}ts)=\{s\}$, by (3), $r\notin \mathcal{R}(w)$.
\\(6) Similar to (5).
\\(7) We use induction on $l(w)$. When $l(w)=0,1,2$,
the lemma is clear. Now assume $l(w)\geq 3$ and the lemma is true for $w'$ if $l(w')<l(w)$. If $w=w_{1}\cdot st=w_{2}\cdot sr$, we have $w=w_{3}\cdot tr$. So $w_{1}s=w_{3}r,\ w_{2}s=w_{3}t$.
We get $w_{3}r=w_{4}\cdot w_{sr},\ w_{3}t=w_{5}\cdot w_{st}$.
Thus $w_{3}=w_{4}\cdot w_{sr}r=w_{5}\cdot w_{st}t$. So $w_{4}\cdot w_{sr}rs=w_{5}\cdot w_{st}ts$, contradict the induction hypothesis.
\\(8) By (1), $\mathcal{L}(tw_{st}w)=\{s\}$, so $r\in \mathcal{L}(rtw_{st}w)$. By (4), $t\notin \mathcal{L}(rtw_{st}w)$.
If $s\in \mathcal{L}(rtw_{st}w)$, we have $rtw_{st}w=w_{sr}\cdot w_{1}$,
then $stw_{st}\cdot w=srw_{sr}\cdot w_{1}$, contradict (7). So $s\notin \mathcal{L}(rtw_{st}w)$.
\\(9) Similar to (8).\hfill $\square$

\medskip

\noindent{\bf Lemma 4.2.}
Let $x,y\in W$.
\\(1) If $w\in W_{st},\ l(w)\geq 4$, $\mathcal{R}(x),\mathcal{L}(y)\subseteq\{r\}$,
then $l(xwy)=l(x)+l(w)+l(y)$.
\\(2) If $w\in W_{sr},\ l(w)\geq 4$, $\mathcal{R}(x),\mathcal{L}(y)\subseteq\{t\}$,
then $l(xwy)=l(x)+l(w)+l(y)$.
\\(3) If $\mathcal{R}(x),\mathcal{L}(y)\subseteq\{s\},\ \mathcal{R}(xt)=\{t\},\ \mathcal{R}(xr)=\{r\},$
then $l(xtry)=l(x)+l(y)+2$.
\\(4) If $\mathcal{R}(x),\mathcal{L}(y)\subseteq \{r\},\ \mathcal{R}(xs)=\{s\}$,
 then $l(xstsy)=l(x)+l(y)+3$.
\\(5) If $\mathcal{R}(x),\mathcal{L}(y)\subseteq \{t\},\ \mathcal{R}(xs)=\{s\}$,
 then $l(xsrsy)=l(x)+l(y)+3$.
\\(6) If $\mathcal{R}(x),\mathcal{L}(y)\subseteq \{r\}$,
then $deg\ (T_{xsts}T_{y})\leq L(r)$.
\\(7) If $\mathcal{R}(x),\mathcal{L}(y)\subseteq\{r\}$, then $l(xtsty)=l(x)+l(y)+3$.

\medskip

\noindent{\bf Proof.}
(1)(2) See [Zhou, 4.4].
\\(3) See the proof of [Zhou, 4.5].
\\(4) By Lemma 4.1(2), $\mathcal{L}(sy)=\{s\}$ or $\mathcal{L}(sy)=\{r,s\}$.
If $\mathcal{L}(sy)=\{s\}$, the lemma is clear.
If $\mathcal{L}(sy)=\{r,s\}$, since $\mathcal{R}(xst)=\{t\}$, the lemma is true by (2).
\\(5) Similar to (4).
\\(6) If $\mathcal{R}(xs)=\{s\}$ or $\mathcal{L}(sy)=\{s\}$, then by (4), $T_{xsts}T_{y}=T_{xstsy}$.
If $\mathcal{R}(xs)=\mathcal{L}(sy)=\{r,s\}$, we have $xs=x_{1}\cdot w_{sr}$, $sy=w_{sr}\cdot y_{1}$,
$\mathcal{R}(x_{1})\subseteq\{t\}$, $\mathcal{L}(y_{1})\subseteq \{t\}$.
By Lemma 4.1(9) and Lemma 4.2(2),
$$
\begin{aligned}
T_{xsts}T_{y}&=T_{x_{1}w_{sr}t}T_{w_{sr}y_{1}}\\
&=(v_{r}-v_{r}^{-1})T_{x_{1}w_{sr}}T_{trw_{sr}y_{1}}+T_{x_{1}w_{sr}r}T_{trw_{sr}y_{1}}\\
&=(v_{r}-v_{r}^{-1})T_{x_{1}w_{sr}trw_{sr}y_{1}}+T_{x_{1}w_{sr}rtrw_{sr}y_{1}}.
\end{aligned}
$$
We get $deg\ (T_{xsts}T_{y})\leq L(r)$.
\\(7) We may assume $y=ry'$, $\mathcal{L}(y')\subseteq\{s\}$.
By Lemma 4.1(5), $\mathcal{R}(xtst)=\{t\}$.
By Lemma 4.1(1), $\mathcal{R}(xts)=\{s\}$.
By Lemma 4.1(4), $t\notin \mathcal{R}(xtsr)$.
By Lemma 4.1(5), $s\notin \mathcal{R}(xtsr)$.
So $\mathcal{R}(xtsr)=\{r\}$.
So by (3), we get $l(xtsty)=l(xtstry')=l(xts)+l(y')+2=l(x)+l(y)+3$.\hfill $\square$

\medskip

\noindent{\bf Lemma 4.3.}
Let $x,y\in W$, $\mathcal{R}(x),\mathcal{L}(y)\subseteq\{s\}$, then $deg\ (T_{xtr}T_{y})\leq L(s)$.

\medskip

\noindent{\bf Proof.}
We have 3 cases.
\\(1) $\mathcal{R}(xt)=\{t\},\ \mathcal{R}(xr)=\{r\},$

By Lemma 4.2(3), $T_{xtr}T_{y}=T_{xtry}$.
\\(2) $\mathcal{R}(xt)=\{s,t\}$.

We assume $xt=x'\cdot w_{st}$, $\mathcal{R}(x')\subseteq\{r\}$. Let $y_{1}=ry$, then $\mathcal{L}(y_{1})=\{r\}$ or $\mathcal{L}(y_{1})=\{s,r\}$. If $\mathcal{L}(y_{1})=\{r\}$, by Lemma 4.2(1), $T_{xtr}T_{y}=T_{x'}T_{w_{st}}T_{_{y_{1}}}=T_{x'w_{st}y_{1}}$. If $\mathcal{L}(y_{1})=\{s,r\}$, we have
$y_{1}=w_{sr}\cdot y'$, $\mathcal{L}(y')\subseteq\{t\}$ for some $y'\in W$. Thus $$T_{xtr}T_{y}=T_{x'w_{st}}T_{w_{sr}y'}=(v_{s}-v_{s}^{-1})T_{x'w_{st}s}T_{w_{sr}y'}+T_{x'w_{st}s}T_{sw_{sr}y'}.$$
By Lemma 4.1(5),(6), $\mathcal{R}(x'w_{st}s)=\{t\},\ \mathcal{L}(sw_{sr}y')=\{r\}$. By Lemma 4.2(2), $T_{x'w_{st}s}T_{w_{sr}y'}=T_{x'w_{st}sw_{sr}y'}$. Since $m_{sr}\geq 5$, by Lemma 4.1(5) and Lemma 4.2(2), $T_{x'w_{st}s}T_{sw_{sr}y'}=T_{x'w_{st}ssw_{sr}y'}$. So we have $deg\ (T_{xtr}T_{y})\leq L(s)$.
\\(3) $\mathcal{R}(xr)=\{s,r\}$.

Similar to (2).\hfill $\square$

\medskip

\noindent{\bf Lemma 4.4.}
Let $x,y\in W,\ \mathcal{R}(x)\subseteq \{t\},\ \mathcal{L}(y)\subseteq\{s\}$,
then we have $deg\ (T_{xw_{sr}}T_{try})\leq L(sr)$.

\medskip

\noindent{\bf Proof.}
We have $T_{xw_{sr}}T_{try}=(v_{r}-v_{r}^{-1})T_{xw_{sr}r}T_{try}+T_{xw_{sr}r}T_{ty}$. Since $\mathcal{R}(xw_{sr}r)=\{s\},\ \mathcal{L}(y)\subseteq\{s\}$, by Lemma 4.3, we get $deg\ ((v_{r}-v_{r}^{-1})T_{xw_{sr}r}T_{try})\leq L(sr)$. We have 2 cases.
\\(1) $\mathcal{L}(ty)=\{t\}$.

By Lemma 4.2(2), $T_{xw_{sr}r}T_{ty}=T_{xw_{sr}rty}$.
\\(2) $\mathcal{L}(ty)=\{s,t\}$.

We have $ty=w_{st}\cdot y'$ for some $y'\in W$.
$$
\begin{aligned}
T_{xw_{sr}r}T_{ty}&=T_{xw_{sr}r}T_{w_{st}y'}\\
&=(v_{s}-v_{s}^{-1})T_{xw_{sr}rs}T_{w_{st}y'}+T_{xw_{sr}rs}T_{sw_{st}y'}\\
&=(v_{s}-v_{s}^{-1})T_{xw_{sr}rsw_{st}y'}+T_{xw_{sr}rsr}T_{rt}T_{tsw_{st}y'}.
\end{aligned}
$$
By Lemma 4.1(1)(2), $\mathcal{R}(xw_{sr}rsr)=\mathcal{L}(tsw_{st}y')=\{s\}$. By Lemma 4.3, $deg\ (T_{xw_{sr}rsr}$ $T_{rt}T_{tsw_{st}y'})\leq L(s)$. So $deg\ (T_{xw_{sr}r}T_{ty})\leq L(s)$ and $deg\ (T_{xw_{sr}}T_{try})\leq L(sr)$.\hfill $\square$

\medskip

\noindent{\bf Lemma 4.5.}
Let $x,y\in W,\ \mathcal{R}(x)\subseteq \{r\},\ \mathcal{L}(y)\subseteq\{s\}$,
then we have
$deg\ (T_{xw_{st}}T_{try})\leq Max$ $\{L(st),L(sr)\}$. Moreover, if $\mathcal{L}(ry)=\{r\}$, then $deg\ (T_{xw_{st}}T_{try})\leq Max\{L(t),L(r)\}$.

\medskip

\noindent{\bf Proof.}
First we have
$$
T_{xw_{st}}T_{try}=(v_{t}-v_{t}^{-1})T_{xw_{st}}T_{ry}+T_{xw_{st}t}T_{ry}.
$$
Since $\mathcal{R}(xw_{st}t)=\{s\},\ \mathcal{L}(y)\subseteq\{s\}$, by Lemma 4.3, we get $deg\ ((v_{t}-v_{t}^{-1})T_{xw_{st}}T_{ry})\leq L(st)$. If $\mathcal{L}(ry)=\{r\}$, by Lemma 4.2(1), $deg\ ((v_{t}-v_{t}^{-1})T_{xw_{st}}T_{ry})\leq L(t)$. Now we consider $T_{xw_{st}t}T_{ry}$, we have 2 cases.
\\(1) $\mathcal{L}(ry)=\{r\}$.

We have 4 cases.
\\\textcircled{\small{1}} $m_{st}\geq 5$.

By Lemma 4.2(1), $T_{xw_{st}t}T_{ry}=T_{xw_{st}try}$.
\\\textcircled{\small{2}} $m_{st}=4,\ \mathcal{R}(xs)=\{s,r\},\ \mathcal{L}(sry)=\{s\}$.

We assume $xs=x_{1}\cdot w_{sr},\ \mathcal{R}(x_{1})\subseteq\{t\}$. Since $\mathcal{L}(sry)=\{s\}$, by Lemma 4.1(3), we have $\mathcal{L}(tsry)=\{t\}$. By Lemma 4.2(2), $T_{xw_{st}t}T_{ry}=T_{xsts}T_{ry}=T_{x_{1}w_{sr}}T_{tsry}=T_{x_{1}w_{sr}tsry}$.
\\\textcircled{\small{3}} $m_{st}=4,\ \mathcal{R}(xs)=\{s,r\},\ \mathcal{L}(sry)=\{s,r\}$.

We assume $xs=x_{1}\cdot w_{sr},\ \mathcal{R}(x_{1})\subseteq\{t\}$. We have $sry=w_{sr}\cdot y_{1},\ \mathcal{L}(y_{1})\subseteq \{t\}$. Since $\mathcal{R}(x_{1}w_{sr}rt)=\{t\}$, $\mathcal{L}(y_{1})\subseteq \{t\}$, by Lemma 4.2(2),
$$
\begin{aligned}
T_{xw_{st}t}T_{ry}&=T_{xsts}T_{ry}\\
&=T_{x_{1}w_{sr}t}T_{w_{sr}y_{1}}\\
&=(v_{r}-v_{r}^{-1})T_{x_{1}w_{sr}rt}T_{w_{sr}y_{1}}+T_{x_{1}w_{sr}rt}T_{rw_{sr}y_{1}}\\
&=(v_{r}-v_{r}^{-1})T_{x_{1}w_{sr}rtw_{sr}y_{1}}+T_{x_{1}w_{sr}rtrw_{sr}y_{1}}.
\end{aligned}
$$
\\\textcircled{\small{4}} $m_{st}=4,\ \mathcal{R}(xs)=\{s\}$.

By Lemma 4.2(4), $T_{xw_{st}t}T_{ry}=T_{xsts}T_{ry}=T_{xstsry}$. So $deg\ (T_{xw_{st}t}T_{ry})\leq L(r)$ if $\mathcal{L}(ry)=\{r\}$.
\\(2) $\mathcal{L}(ry)=\{r,s\}$.

We assume $ry=w_{sr}\cdot y',\ \mathcal{L}(y')\subseteq\{t\}$. So
$$
T_{xw_{st}t}T_{ry}=T_{xw_{st}t}T_{w_{sr}y'}
=(v_{s}-v_{s}^{-1})T_{xw_{st}ts}T_{w_{sr}y'}+T_{xw_{st}ts}T_{sw_{sr}y'}.
$$

We have 2 cases.
\\\textcircled{\small{1}} $\mathcal{R}(xw_{st}ts)=\{t\}$.

By Lemma 4.2(2), $T_{xw_{st}ts}T_{w_{sr}y'}=T_{xw_{st}tsw_{sr}y'}$. On the other hand, $T_{xw_{st}ts}T_{sw_{sr}y'}$ $=T_{xw_{st}tst}T_{tr}T_{rsw_{sr}y'}$. Since $\mathcal{R}(xw_{st}tst)=\mathcal{L}(rsw_{sr}y')=\{s\}$, by Lemma 4.3, we get $deg\ (T_{xw_{st}ts}T_{sw_{sr}y'})\leq L(s)$ and $deg\ (T_{xw_{st}t}T_{ry})$ $\leq L(s)$.
\\\textcircled{\small{2}} $\mathcal{R}(xw_{st}ts)=\{r,t\}$.

Then $m_{st}=4$ and $\mathcal{R}(xs)=\{r,s\}$. We have $xs=x'\cdot w_{sr},\ \mathcal{R}(x')\subseteq\{t\}$. So by Lemma 4.2(2),
$$
\begin{aligned}
T_{xw_{st}t}T_{ry}&=T_{xw_{st}t}T_{w_{sr}y'}\\
&=(v_{s}-v_{s}^{-1})T_{xst}T_{w_{sr}y'}+T_{xst}T_{sw_{sr}y'}\\
&=(v_{s}-v_{s}^{-1})T_{x'w_{sr}t}T_{w_{sr}y'}+T_{x'w_{sr}t}T_{sw_{sr}y'}\\
&=(v_{s}-v_{s}^{-1})(v_{r}-v_{r}^{-1})T_{x'w_{sr}rtw_{sr}y'}
+(v_{r}-v_{r}^{-1})T_{x'w_{sr}rtrw_{sr}y'}
\\&\ \ \ \ +(v_{r}-v_{r}^{-1})T_{x'w_{sr}rtsw_{sr}y'}
+T_{x'w_{sr}rtrsw_{sr}y'}.
\end{aligned}
$$
Thus, $deg\ (T_{xw_{st}t}T_{ry})\leq L(sr)$. We get $deg\ (T_{xw_{st}}T_{try})\leq Max\{L(st),$ $L(sr)\}$. Moreover, $deg\ (T_{xw_{st}}T_{try})\leq Max\{L(t),L(r)\}$ if $\mathcal{L}(ry)=\{r\}$.\hfill $\square$

\medskip

\noindent{\bf Lemma 4.6.}
Let $c\in W_{st}$, $l(c)\leq m_{st}-2$ or $c=sw_{st}$, then we have $deg\ f_{w_{st},c,st}\leq L(t)$, $deg\ f_{w_{st},c,ts}\leq L(t)$, $deg\ f_{w_{st},c,tst}\leq 2L(t)$, $deg\ f_{w_{st},c,sts}\leq L(st)$.

\medskip

\noindent{\bf Proof.}
We have 3 cases.
\\(1) $l(c)\leq m_{st}-3$.

We have $f_{w_{st},c,st}=f_{w_{st},c,ts}=0$, $f_{w_{st},c,tst}=0$ or $1$, $f_{w_{st},c,sts}=0$ or $1$ in this case.
\\(2) $l(c)=m_{st}-2$.

We have $f_{w_{st},c,st}=0$ or $1$, $f_{w_{st},c,ts}=0$ or $1$, $f_{w_{st},c,tst}=v_{t}-v_{t}^{-1}$, $f_{w_{st},c,sts}=v_{s}-v_{s}^{-1}$ in this case.
\\(3) $c=sw_{st}$.

We have $$(v_{s}-v_{s}^{-1})T_{w_{st}}T_{sw_{st}}=T_{w_{st}}T_{w_{st}}-T_{w_{st}s}T_{sw_{st}}.$$
Thus, $(v_{s}-v_{s}^{-1})f_{w_{st},sw_{st},q}=f_{w_{st},w_{st},q}-f_{w_{st}s,sw_{st},q}$ for any $q \in W_{st}$. By Lemma 1.4(2), we have $deg\ f_{w_{st},sw_{st},st}\leq L(t)$, $deg\ f_{w_{st},sw_{st},ts}\leq L(t)$, $deg\ f_{w_{st},sw_{st},tst}\leq 2L(t)$, $deg\ f_{w_{st},sw_{st},sts}\leq L(st)$.\hfill $\square$

\medskip

\noindent{\bf Proposition 4.7.}
We have $deg\ (T_{x}T_{y})\leq N$ for all $x,y \in W$ in case 2.

\noindent{\bf Proof.}
Now we prove Theorem 2.1 for this case. We use induction on $l(y)$. When $l(y)=0,1,2$, the proposition is clear by Lemma 1.4(2). Now assume $l(y)\geq 3$
and the proposition is true for $y'$ if $l(y')<l(y)$. We assume $x=x'\cdot w$, $y=u\cdot y'$, $w,u\in W_{sr}$, $\mathcal{R}(x'),\mathcal{L}(y')\subseteq \{t\}$. Then
$$
\begin{aligned}
T_{x}T_{y}&=T_{x'}(T_{w}T_{u})T_{y'}\\
&=T_{x'}(\sum_{p\in W_{sr}}f_{w,u,p}T_{p})T_{y'}\\
&=\sum_{p\in W_{sr}}f_{w,u,p}T_{x'p}T_{y'}.
\end{aligned}
$$

We will prove $deg\ (f_{w,u,p}T_{x'p}T_{y'})\leq N$ for all $p\in W_{sr}$. By Lemma 1.4(2), $deg\ f_{w,u,p}\leq L(p)$. If $l(p)\geq 4$, by Lemma 4.2(2), $T_{x'p}T_{y'}=T_{x'py'}$, so $deg\ (f_{w,u,p}T_{x'p}T_{y'})$ $\leq N$.
If $l(p)\leq 3$, $l(x')\leq 1$ or $l(y')\leq 1$, it is easy to see $deg\ (f_{w,u,p}T_{x'p}T_{y'})\leq N$. From now on we assume $l(p)\leq 3$, $x'=x_{1}\cdot st$, $y'=ts\cdot y_{1}$, $\mathcal{R}(x_{1}s)=\mathcal{L}(sy_{1})=\{s\}$. We have 7 cases.
\\(1) $p=e$.

By Lemma 1.4(1), $f_{w,u,e}\in \mathbb{Z}$, so we only need to prove $deg\ (T_{x'}T_{y'})\leq N$. We have $$T_{x'}T_{y'}=T_{x_{1}st}T_{tsy_{1}}=(v_{t}-v_{t}^{-1})T_{x_{1}st}T_{sy_{1}}+T_{x_{1}s}T_{sy_{1}}.$$ Since $l(sy_{1})<l(y)$, by induction hypothesis, we have $deg\ (T_{x_{1}s}T_{sy_{1}})\leq N$, so we only need to prove $deg\ (T_{x_{1}st}T_{sy_{1}})\leq N-L(t)$. We assume $x_{1}st=x_{2}\cdot a,\ sy_{1}=b\cdot y_{2},\ a,b\in W_{st},\ \mathcal{R}(x_{2}),\ \mathcal{L}(y_{2})\subseteq\{r\}$. Then we have $3\leq l(a)+l(b) \leq 2m_{st}-3$.

If $l(a)+l(b)=3$, by Lemma 4.2(4), $T_{x_{1}st}T_{sy_{1}}=T_{x_{1}stsy_{1}}$.

If $4\leq l(a)+l(b)\leq m_{st}$, by Lemma 4.2(1), $T_{x_{1}st}T_{sy_{1}}=T_{x_{1}stsy_{1}}$.

If $m_{st}+1\leq l(a)+l(b)\leq 2m_{st}-3$, take $c\in W_{st}$ such that $ab=w_{st}c$ and $1\leq l(c)\leq m_{st}-3$, then
$$
\begin{aligned}
T_{x_{1}st}T_{sy_{1}}&=T_{x_{2}a}T_{by_{2}}\\
&=T_{x_{2}w_{st}}T_{cy_{2}}\\
&=T_{x_{2}}(\sum_{q\in W_{st},l(q)\geq 3}f_{w_{st},c,q}T_{q})T_{y_{2}}\\
&=f_{w_{st},c,sts}T_{x_{2}sts}T_{y_{2}}+f_{w_{st},c,tst}T_{x_{2}tsty_{2}}+\sum_{q\in W_{st},l(q)\geq 4}f_{w_{st},c,q}T_{x_{2}qy_{2}}.
\end{aligned}
$$
It is clear that $f_{w_{st},c,sts}, f_{w_{st},c,tst}\in \{0,1\}$. By Lemma 4.2(6), we have $deg\ (T_{x_{2}sts}T_{y_{2}})\leq L(r)$. By Lemma 1.4(2), $deg\ f_{w_{st},c,q}\leq L(c)\leq N-L(t)$ for any $q\in W_{st}$. So we get $deg\ (T_{x_{1}st}T_{sy_{1}})\leq N-L(t)$.
\\(2) $p=s$.

Since $deg\ f_{w,u,s}\leq L(s)$, we only need to prove $deg\ (T_{x's}T_{y'})\leq N-L(s)$. We assume $x's=x_{1}sts=x_{2}\cdot a,\ y'=tsy_{1}=b\cdot y_{2}$,
$a,b\in W_{st},\ \mathcal{R}(x_{2}),\mathcal{L}(y_{2})\subseteq\{r\}$,
then we have $5\leq l(a)+l(b)\leq 2m_{st}-1$.

If $5\leq l(a)+l(b)\leq m_{st}$, by Lemma 4.2(1), $T_{x's}T_{y'}=T_{x_{2}a}T_{by_{2}}=T_{x_{2}aby_{2}}$.

If $m_{st}+1\leq l(a)+l(b)\leq 2m_{st}-1$, then $T_{x's}T_{y'}=T_{x_{2}a}T_{by_{2}}=T_{x_{2}w_{st}}T_{cy_{2}}$, $c\in W_{st}$ such that
$ab=w_{st}c$ and $1\leq l(c)\leq m_{st}-1$,
$$
T_{x_{2}w_{st}}T_{cy_{2}}=T_{x_{2}}(\sum_{q\in W_{st}\setminus \{e\}}f_{w_{st},c,q}T_{q})T_{y_{2}}=\sum_{q\in W_{st}\setminus \{e\}}f_{w_{st},c,q}T_{x_{2}q}T_{y_{2}}.
$$
We will prove $deg\ (f_{w_{st},c,q}T_{x_{2}q}T_{y_{2}})\leq N-L(s)$ for any $q\in W_{st}\setminus \{e\}$.
\\\textcircled{\small{1}} $q=s$.

If $f_{w_{st},c,q}\neq 0$, then $a=w_{st}$, $b=c=sw_{st}$, $f_{w_{st},c,q}=1$, $m_{st}$ is even. We have $(v_{s}-v_{s}^{-1})T_{x_{2}s}T_{y_{2}}=T_{x_{2}s}T_{sy_{2}}-T_{x_{2}}T_{y_{2}}$. Since $l(sy_{2})<l(y)$, $l(y_{2})<l(y)$, by induction hypothesis, we get $deg\ (T_{x_{2}s}T_{sy_{2}})\leq N$ and $deg\ (T_{x_{2}}T_{y_{2}})\leq N$. So $deg\ (T_{x_{2}s}T_{y_{2}})\leq N-L(s)$.
\\\textcircled{\small{2}} $q=t$.

If $f_{w_{st},c,q}\neq 0$, then $a=w_{st}$, $b=c=sw_{st}$, $f_{w_{st},c,q}=1$, $m_{st}$ is odd. We have $(v_{s}-v_{s}^{-1})T_{x_{2}t}T_{y_{2}}=(v_{t}-v_{t}^{-1})T_{x_{2}t}T_{y_{2}}=T_{x_{2}t}T_{ty_{2}}-T_{x_{2}}T_{y_{2}}$. Since $l(ty_{2})<l(y)$, $l(y_{2})<l(y)$, by induction hypothesis, we get $deg\ (T_{x_{2}t}T_{ty_{2}})\leq N$ and $deg\ (T_{x_{2}}T_{y_{2}})\leq N$. So $deg\ (T_{x_{2}t}T_{y_{2}})\leq N-L(s)$.
\\\textcircled{\small{3}} $q=st$, $l(y_{2})\leq 1$.

By Lemma 4.6, we have $deg\ f_{w_{st},c,q}\leq L(t)$. Since $deg\ (T_{x_{2}q}T_{y_{2}})\leq L(r)$, we get $deg\ (f_{w_{st},c,q}T_{x_{2}q}T_{y_{2}})\leq L(rt)$.
\\\textcircled{\small{4}} $q=st$, $l(y_{2})\geq 2$.

By Lemma 4.6, we have $deg\ f_{w_{st},c,q}\leq L(t)$. Assume $y_{2}=rsy_{3}$, $\mathcal{L}(sy_{3})=\{s\}$, by Lemma 4.3, Lemma 4.4, we get $deg\ (T_{x_{2}q}T_{y_{2}})=deg\ (T_{x_{2}st}T_{rsy_{3}})\leq L(sr)$, so $deg\ (f_{w_{st},c,q}T_{x_{2}q}T_{y_{2}})\leq L(srt)$.
\\\textcircled{\small{5}} $q=ts$, $l(x_{2})\leq 1$.

Similar to \textcircled{\small{3}}.
\\\textcircled{\small{6}} $q=ts$, $l(x_{2})\geq 2$.

Similar to \textcircled{\small{4}}.
\\\textcircled{\small{7}} $q=tst$.

By Lemma 4.6, we have $deg\ f_{w_{st},c,q}\leq 2L(t)$. By Lemma 4.2(7), $T_{x_{2}q}T_{y_{2}}=T_{x_{2}tst}T_{y_{2}}=T_{x_{2}tsty_{2}}$. So $deg\ (f_{w_{st},c,q}T_{x_{2}q}T_{y_{2}})\leq 2L(t)$.
\\\textcircled{\small{8}} $q=sts$.

By Lemma 4.6, we have $deg\ f_{w_{st},c,q}\leq L(st)$. By Lemma 4.2(6), we have $deg\ (T_{x_{2}q}T_{y_{2}})\leq L(r)$. So $deg\ (f_{w_{st},c,q}T_{x_{2}q}T_{y_{2}})\leq L(str)$.
\\\textcircled{\small{9}} $l(q)\geq 4$.

By Lemma 4.2(1), we have $deg\ (f_{w_{st},c,q}T_{x_{2}q}T_{y_{2}})=deg\ f_{w_{st},c,q}\leq L(c)\leq N-L(s)$.
\\(3) $p=r$.

We have 2 cases.
\\\textcircled{\small{1}} $u=e$.

We have
$$T_{x'r}T_{y'}=(v_{t}-v_{t}^{-1})T_{x_{1}str}T_{sy_{1}}+T_{x_{1}sr}T_{sy_{1}}.$$
By Lemma 4.3, we have $deg\ ((v_{t}-v_{t}^{-1})T_{x_{1}str}T_{sy_{1}})\leq L(st)$.
Since $l(sy_{1})<l(y)$, by induction hypothesis, we have $deg\ (T_{x_{1}sr}T_{sy_{1}})\leq N$. So $deg\ (T_{x'r}T_{y'})\leq N$ and $deg\ (f_{w,e,r}T_{x'r}T_{y'})\leq N$.
\\\textcircled{\small{2}} $u\neq e$.

We have
$$
\begin{aligned}
&\ \ \ \ (v_{r}-v_{r}^{-1})T_{x'r}T_{y'}\\
&=(v_{r}-v_{r}^{-1})T_{x_{1}str}T_{tsy_{1}}\\
&=(v_{r}-v_{r}^{-1})(v_{t}-v_{t}^{-1})T_{x_{1}str}T_{sy_{1}}+
(v_{r}-v_{r}^{-1})T_{x_{1}sr}T_{sy_{1}}\\
&=(v_{r}-v_{r}^{-1})(v_{t}-v_{t}^{-1})T_{x_{1}str}T_{sy_{1}}+T_{x_{1}sr}T_{rsy_{1}}-T_{x_{1}s}T_{sy_{1}}.
\end{aligned}
$$
By Lemma 4.3, we know $deg\ ((v_{r}-v_{r}^{-1})(v_{t}-v_{t}^{-1})T_{x_{1}str}T_{sy_{1}})\leq L(rst)$. Since $l(rsy_{1})<l(y)$, $l(sy_{1})<l(y)$, by induction hypothesis, we get $deg\ (T_{x_{1}sr}T_{rsy_{1}})\leq N$,
$deg\ (T_{x_{1}s}T_{sy_{1}})\leq N$. Thus, we have $deg\ (T_{x'r}T_{y'})\leq N-L(r)$. Since $deg\ f_{w,u,r}\leq L(r)$, we get $deg\ (f_{w,u,r}T_{x'r}T_{y'})\leq N$.
\\(4) $p=sr$.

First we know $deg\ f_{w,u,sr}\leq L(sr)$. We have $T_{x'sr}T_{y'}=T_{x's}T_{rt}T_{sy_{1}}$. If $\mathcal{R}(x's)=\{s,t\}$, by Lemma 4.5, $deg\ (T_{x'sr}T_{y'})\leq Max\{L(st),L(sr)\}$, so $deg\ (f_{w,u,sr}$ $T_{x'sr}T_{y'})\leq Max\{L(srst),L(srsr)\}$. If $\mathcal{R}(x's)=\{s\}$, by Lemma 4.3, $deg\ (T_{x'sr}T_{y'})\leq L(s)$, so $deg\ (f_{w,u,sr}T_{x'sr}T_{y'})\leq L(srs)$.
\\(5) $p=rs$.

Similar to (4).
\\(6) $p=rsr$.

We have $T_{x'rsr}T_{y'}=T_{x'rs}T_{rt}T_{sy_{1}}$. By Lemma 4.1(2), we get $\mathcal{R}(x'rs)=\mathcal{L}(sy_{1})=\{s\}$, so by Lemma 4.3, $deg\ (T_{x'rsr}T_{y'})\leq L(s)$. Since $deg\ f_{w,u,rsr}$ $\leq L(rsr)$, we get $deg\ (f_{w,u,rsr}T_{x'rsr}T_{y'})\leq L(rsrs)$.
\\(7) $p=srs$.

First we have $deg\ f_{w,u,srs}\leq L(srs)$. If $\mathcal{R}(x's)=\{s\}$ or $\mathcal{L}(sy')=\{s\}$, by Lemma 4.2(5), $T_{x'srs}T_{y'}=T_{x'srsy'}$,
so $deg\ (f_{w,u,srs}T_{x'srs}T_{y'})\leq L(srs)$. If $\mathcal{R}(x's)=\{s,t\}$ and $\mathcal{L}(sy')=\{s,t\}$, we assume $x's=x''\cdot w_{st},\ sy'=w_{st}\cdot y'',\ \mathcal{R}(x''),\mathcal{L}(y'')\subseteq \{r\}$. Then $T_{x'srs}T_{y'}=T_{x''w_{st}r}T_{w_{st}y''}=T_{x''w_{st}}T_{tr}T_{tw_{st}y''}$. Since $\mathcal{R}(x'')\subseteq\{r\}$, $\mathcal{L}(tw_{st}y'')=\{s\}$, $\mathcal{L}(rtw_{st}y'')=\{r\}$, by Lemma 4.5, we get $deg\ (T_{x''w_{st}}T_{tr}$ $T_{tw_{st}y''})\leq Max\{L(t),L(r)\}$.
So $deg\ (f_{w,u,srs}T_{x'srs}T_{y'})\leq Max\{L(srst),$ $L(srsr)\}$.\hfill $\square$

\medskip

\section{Case 3}

\medskip

In this section, we prove Theorem 2.1 for the case $\infty>m_{sr}\geq7$, $m_{st}=3$. In this case, $N=L(w_{sr})>Max\{L(w_{st}),L(rt)\}$.

\medskip

\noindent{\bf Lemma 5.1.}
Let $w\in W$.
\\(1) There is no $w_{1},w_{2}\in W$ such that $w=w_{1}\cdot st=w_{2}\cdot sr$.
\\(2) If $w=w_{1}\cdot srs$, then $t\notin \mathcal{R}(w)$.
\\(3) If $w=w_{1}\cdot srsr$, then $t\notin \mathcal{R}(w)$.
\\(4) If $w=w_{1}\cdot ts$, then $r\notin \mathcal{R}(w)$.
\\(5) If $w=w_{1}\cdot tsr$, then $s\notin \mathcal{R}(w)$.

\medskip

\noindent{\bf Proof.}
(1) We use induction on $l(w)$. It is easy to check this lemma when $l(w)\leq 5$. Now assume $l(w)\geq 6$ and the lemma is true for $w'$ if $l(w')<l(w)$.
We assume $w=w_{3}\cdot rt$. So $w_{1}s=w_{3}r$, $w_{2}s=w_{3}t$. We have $w_{1}s=w_{3}r=w_{4}\cdot w_{sr}$,
$w_{2}s=w_{3}t=w_{5}\cdot w_{st}$. Since $m_{sr}\geq 7$, we assume $w_{3}=w_{5}\cdot w_{st}t=w_{5}\cdot ts=w_{4}\cdot w_{sr}r=w_{4}'\cdot srsrs$. Since $\mathcal{R}(w_{5})\subseteq\{r\}$ and $l(w_{5})\geq 2$, we assume $w_{5}=w_{5}'\cdot sr$. Thus $w_{4}'\cdot srs=w_{5}'\cdot st$. We get $w_{4}'\cdot srs=w_{6}\cdot sts$,
so $w_{4}'\cdot sr=w_{6}\cdot st$, contradict the induction hypothesis.
\\(2) If $w=w_{1}\cdot srs=w_{2}\cdot t$, we have $w=w_{1}\cdot srs=w_{3}\cdot sts$, so $w_{1}\cdot sr=w_{3}\cdot st$, contradict (1).
\\(3) If $w=w_{1}\cdot srsr=w_{2}\cdot t$, we have $w=w_{1}\cdot srsr=w_{3}\cdot tr$, so $w_{1}\cdot srs=w_{3}\cdot t$, contradict (2).
\\(4) If $w=w_{1}\cdot ts=w_{2}\cdot r$, we have $w=w_{1}\cdot ts=w_{3}\cdot w_{sr}$, so $w_{1}\cdot t=w_{3}\cdot w_{sr}s$, contradict (3).
\\(5) If $s\in \mathcal{R}(w)$, then we have $w=w_{2}\cdot srsrsr$ for some $w_{2}\in W$, so $w_{1}\cdot t=w_{2}\cdot srsr$, contradict (3).\hfill $\square$

\medskip

According to [Zhou, 3.5], we have
\\{\bf Lemma 5.2.}
Let $x,y\in W$, $\mathcal{R}(x),\mathcal{L}(y)\subseteq\{t\}$, $w\in W_{sr}$, $l(w)\geq 6$ or $w=srsrs$,
then $l(xwy)=l(x)+l(w)+l(y)$, $\mathcal{R}(xwy)=\mathcal{R}(wy)$, $\mathcal{L}(xwy)=\mathcal{L}(xw)$.

\medskip

\noindent{\bf Lemma 5.3.}
Let $x,y\in W$, $\mathcal{R}(x),\mathcal{L}(y)\subseteq \{r\}$.
\\(1) If $\mathcal{R}(xs)=\{s\}$ or $\mathcal{L}(sy)=\{s\}$, then $T_{xsts}T_{y}=T_{xstsy}$.
\\(2) If $\mathcal{R}(xs)=\mathcal{L}(sy)=\{r,s\}$, then $deg\ (T_{xsts}T_{y})=L(r)$.

\medskip

\noindent{\bf Proof.}
(1) See the proof of [Zhou, 3.6].
\\(2) We assume $x=x'\cdot w_{sr}s$, $y=sw_{sr}\cdot y'$ for some $x',y'\in W$ with $\mathcal{R}(x'),\mathcal{L}(y')\subseteq\{t\}$. By Lemma 5.2, we have
$$
\begin{aligned}
T_{xsts}T_{y}&=T_{x'\cdot w_{rs}}T_{t}T_{w_{rs}\cdot y'}\\
&=\xi_{r}T_{x'\cdot w_{rs}\cdot t\cdot rw_{rs}\cdot y'}+T_{x'\cdot w_{rs}r\cdot t\cdot rw_{rs}\cdot y'}.
\end{aligned}
$$
\hfill $\square$

\medskip

\noindent{\bf Lemma 5.4.}
Let $x,y\in W$, $\mathcal{R}(x),\mathcal{L}(y)\subseteq \{s\}$.
\\(1) If $\mathcal{R}(xr)\neq \{s,r\}$, $\mathcal{R}(xt)\neq \{s,t\}$, $\mathcal{R}(xrs)\neq \{s,r\}$, then $T_{xtr}T_{y}=T_{xtry}$.
\\(2) If $\mathcal{R}(xr)=\{s,r\}$, then $deg\ (T_{xtr}T_{y})\leq L(sr)$.
\\(3) If $\mathcal{R}(xt)=\{s,t\}$, then $deg\ (T_{xtr}T_{y})\leq L(sr)$.
\\(4) If $\mathcal{R}(xrs)=\{s,r\}$, then $deg\ (T_{xtr}T_{y})\leq L(r)$.

\medskip

\noindent{\bf Proof.}
(1) See the proof of [Zhou, 3.7].
\\(2) We assume $xr=x'\cdot w_{sr}$ for some $x'\in W$, then $xtr=x'\cdot w_{sr}\cdot t$, so $T_{xtr}T_{y}=T_{x'w_{sr}}T_{ty}$.
We have $ty=u\cdot y'$, $u\in W_{sr}$, $\mathcal{L}(y')\subseteq\{t\}$. By Lemma 5.1(2), $u=e$ or $u=s$ or $u=sr$.
If $u=e$, then $T_{xtr}T_{y}=T_{x'w_{sr}}T_{y'}=T_{x'w_{sr}y'}$.
If $u=s$, then $T_{xtr}T_{y}=T_{x'w_{sr}}T_{sy'}=(v_{s}-v_{s}^{-1})T_{x'w_{sr}y'}+T_{x'w_{sr}sy'}$.
If $u=sr$, then
$$
\begin{aligned}
T_{xtr}T_{y}&=T_{x'w_{sr}}T_{sry'}\\
&=(v_{s}-v_{s}^{-1})(v_{r}-v_{r}^{-1})T_{x'w_{sr}y'}+(v_{s}-v_{s}^{-1})T_{x'w_{sr}ry'}\\
&\ \ \ \ +(v_{r}-v_{r}^{-1})T_{x'w_{sr}sy'}+T_{x'w_{sr}sry'}.
\end{aligned}
$$
\\(3) We assume $xt=x'\cdot w_{st}$ for some $x'\in W$, then $xtr=x'\cdot w_{st}\cdot r$, so $T_{xtr}T_{y}=T_{x'w_{st}}T_{ry}$.
If $\mathcal{L}(ry)=\{r\}$, then by Lemma 5.3, $deg\ (T_{xtr}T_{y})\leq L(r)$.
If $\mathcal{L}(ry)=\{s,r\}$, then we have $ry=w_{sr}\cdot u$ for some $u\in W$, so $T_{x'w_{st}}T_{ry}=T_{x'w_{st}}T_{w_{sr}u}
=(v_{s}-v_{s}^{-1})T_{x'w_{st}}T_{sw_{sr}u}+T_{x'w_{st}s}T_{sw_{sr}u}$. Since $\mathcal{R}(x')=\mathcal{L}(sw_{sr}u)=\{r\}$, by Lemma 5.3, $deg\ ((v_{s}-v_{s}^{-1})T_{x'w_{st}}T_{sw_{sr}u})\leq L(sr)$. If $\mathcal{R}(x'w_{st}s)=\{t\}$, then by Lemma 5.2,
$T_{x'w_{st}s}T_{sw_{sr}u}=T_{x'stsw_{sr}u}$. If $\mathcal{R}(x'w_{st}s)=\{r,t\}$, then we get $x's=x''\cdot w_{sr}$ for some $x''\in W$. So
$$
\begin{aligned}
T_{x'w_{st}s}T_{sw_{sr}u}&=T_{x''w_{sr}t}T_{sw_{sr}u}\\
&=(v_{r}-v_{r}^{-1})T_{x''w_{sr}t}T_{rsw_{sr}u}+T_{x''w_{sr}rt}T_{rsw_{sr}u}\\
&=(v_{r}-v_{r}^{-1})T_{x''w_{sr}trsw_{sr}u}+T_{x''w_{sr}rtrsw_{sr}u}.
\end{aligned}
$$
\\(4) We assume $xrs=x'\cdot w_{sr}$ for some $x'\in W$, then $xtr=x'\cdot w_{sr}s\cdot t$, so $T_{xtr}T_{y}=T_{x'w_{sr}s}T_{ty}$.
We have $ty=u\cdot y'$, $u\in W_{sr}$, $\mathcal{L}(y')\subseteq\{t\}$. By Lemma 5.1(2), $u=e$ or $u=s$ or $u=sr$.
If $u=e$, then $T_{xtr}T_{y}=T_{x'w_{sr}s}T_{y'}=T_{x'w_{sr}sy'}$.
If $u=s$, then $T_{xtr}T_{y}=T_{x'w_{sr}s}T_{sy'}=T_{x'w_{sr}y'}$.
If $u=sr$, then $T_{xtr}T_{y}=T_{x'w_{sr}s}T_{sry'}=(v_{r}-v_{r}^{-1})T_{x'w_{sr}y'}+T_{x'w_{sr}ry'}$.\hfill $\square$

\medskip

\noindent{\bf Lemma 5.5.}
Let $x,y\in W$, $\mathcal{R}(x)\subseteq \{r\}$, $\mathcal{L}(y)\subseteq \{t\}$, then we have $deg\ (T_{xw_{st}}T_{w_{sr}y})\leq L(sr)$.

\medskip

\noindent{\bf Proof.}
See the proof of Lemma 5.4(3).\hfill $\square$

\medskip

\noindent{\bf Lemma 5.6.}
Let $x,y\in W$, $\mathcal{R}(x),\mathcal{L}(y)\subseteq \{t\}$, then we have $deg\ (T_{xrsrsr}$ $T_{y})\leq L(s)$.

\medskip

\noindent{\bf Proof.}
We suppose $l(xrsrsry)<l(x)+l(y)+5$. Then we have $l(x),l(y)\geq 1$, so we may assume $x=x'\cdot t$ and $y=t\cdot y'$ for some $x',y'\in W$ with $\mathcal{R}(x'),\mathcal{L}(y')\subseteq\{s\}$, then $T_{x}T_{rsrsr}T_{y}=T_{x'}T_{rt}T_{srsrty'}$.
Since $\mathcal{L}(rsrsrty')=\{r\}$, $\mathcal{L}(srsrsrty')=\{s\}$, by Lemma 5.4, we must have $\mathcal{L}(tsrsrty')=\{s,t\}$. Thus we get $\mathcal{L}(ry')=\{r,s\}$. Similarly, we can prove $\mathcal{R}(x'r)=\{r,s\}$.
Now we assume $x=x''\cdot w_{sr}r\cdot t$, $y=t\cdot rw_{sr}\cdot y''$ for some $x'',y''\in W$ with $\mathcal{R}(x''),\mathcal{L}(y'')\subseteq\{t\}$.
By Lemma 5.2, we have $\mathcal{L}(tsrst\cdot sw_{sr}\cdot y'')=\mathcal{L}(tsrst\cdot sw_{sr})=\{t\}$, thus
$$
\begin{aligned}
T_{x}T_{rsrsr}T_{y}&=T_{x''\cdot w_{sr}r\cdot t}T_{rsrsr}T_{t\cdot rw_{sr}\cdot y''}\\
&=T_{x''\cdot w_{sr}s\cdot tst}T_{r}T_{tst\cdot sw_{sr}\cdot y''}\\
&=\xi_{t}T_{x''\cdot w_{sr}\cdot tsrst\cdot sw_{sr}\cdot y''}+T_{x''\cdot w_{sr}s\cdot tsrst\cdot sw_{sr}\cdot y''}.
\end{aligned}
$$
\hfill $\square$

\medskip

\noindent{\bf Lemma 5.7.}
Let $x,y\in W$.
\\(1) If $\mathcal{R}(x)\subseteq\{s\}$, $\mathcal{L}(y)\subseteq\{t\}$, then $deg\ (T_{xtr}T_{w_{sr}y})\leq L(rsr)$.
\\(2) If $\mathcal{R}(x)\subseteq\{t\}$, $\mathcal{L}(y)\subseteq\{s\}$, then $deg\ (T_{xw_{sr}}T_{try})\leq L(rsr)$.

\medskip

\noindent{\bf Proof.}
(1) We have $$T_{xtr}T_{w_{sr}y}=(v_{r}-v_{r}^{-1})T_{xtr}T_{rw_{sr}y}+T_{xt}T_{rw_{sr}y}.$$
Since $\mathcal{R}(x)\subseteq\{s\}$, $\mathcal{L}(rw_{sr}y)\subseteq\{s\}$,
by Lemma 5.4, we have $deg\ ((v_{r}-v_{r}^{-1})T_{xtr}T_{rw_{sr}y})$ $\leq L(rsr)$. If $\mathcal{R}(xt)=\{t\}$,
by Lemma 5.2, $T_{xt}T_{rw_{sr}y}=T_{xtrw_{sr}y}$. If $\mathcal{R}(xt)=\{s,t\}$, we have $xt=x'\cdot w_{st}$ for some $x'\in W$,
then $T_{xt}T_{rw_{sr}y}=(v_{s}-v_{s}^{-1})T_{x'sts}T_{srw_{sr}y}+T_{x'st}T_{srw_{sr}y}$. By Lemma 5.3, we have $deg\ ((v_{s}-v_{s}^{-1})T_{x'sts}T_{srw_{sr}y})\leq L(rs)$. If $\mathcal{R}(x'st)=\{t\}$,
by Lemma 5.2 and Lemma 5.6, $deg\ (T_{x'st}T_{srw_{sr}y})\leq L(s)$. If $\mathcal{R}(x'st)=\{r,t\}$, we get $x'st=x''\cdot w_{sr}\cdot t$,
thus
$$
\begin{aligned}
T_{x'st}T_{srw_{sr}y}&=T_{x''w_{sr}t}T_{srw_{sr}y}\\
&=(v_{r}-v_{r}^{-1})T_{x''w_{sr}rt}T_{srw_{sr}y}+T_{x''w_{sr}rt}T_{rsrw_{sr}y}\\
&=(v_{r}-v_{r}^{-1})T_{x''w_{sr}rtsrw_{sr}y}+T_{x''w_{sr}rtrsrw_{sr}y}.
\end{aligned}
$$
\\(2) Apply the $\mathbb{Z}[v,v^{-1}]-$algebra antiautomorphism $T_{w}\longrightarrow T_{w^{-1}}$ of $\mathcal{H}$ to (1).\hfill $\square$

\medskip

\noindent{\bf Lemma 5.8.}
Let $x,y\in W$.
\\(1) If $\mathcal{R}(x)\subseteq\{s\}$, $\mathcal{L}(y)\subseteq \{r\}$, then $deg\ (T_{xtr}T_{w_{st}y})\leq L(srsr)$.
\\(2) If $\mathcal{R}(x)\subseteq\{r\}$, $\mathcal{L}(y)\subseteq \{s\}$, then $deg\ (T_{xw_{st}}T_{try})\leq L(srsr)$.

\medskip

\noindent{\bf Proof.}
(1) We have $$T_{xtr}T_{tsty}=(v_{s}-v_{s}^{-1})T_{xtr}T_{sty}+T_{xr}T_{sty}.$$
Since $\mathcal{L}(sty)=\{s\}$, by Lemma 5.4, we have
$deg\ ((v_{s}-v_{s}^{-1})T_{xtr}T_{sty})\leq L(srs)$. Now we consider $T_{xr}T_{sty}$. If $y=e$, the lemma is clear, so we may assume $y=r\cdot y_{1}$ for
some $y_{1}\in W$, $\mathcal{L}(y_{1})\subseteq\{s\}$. We have 2 cases.
\\\textcircled{\small{1}} $\mathcal{R}(xr)=\{r\}$.

If $\mathcal{R}(xrs)=\{s\}$, then $T_{xr}T_{sty}=T_{xrstr}T_{y_{1}}$. By Lemma 5.4, we have $deg\ (T_{xr}T_{sty})$ $\leq L(sr)$.
If $\mathcal{R}(xrs)=\{s,r\}$, then we assume $xrs=x'\cdot w_{sr}$ for some $x'\in W$. So $T_{xr}T_{sty}=T_{x'w_{sr}}T_{try_{1}}$. By Lemma 5.7(2),
$deg\ (T_{xr}T_{sty})\leq L(rsr)$.
\\\textcircled{\small{2}} $\mathcal{R}(xr)=\{s,r\}$.

We assume $xr=x'\cdot w_{sr}$ for some $x'\in W$, thus
$$
\begin{aligned}
T_{xr}T_{sty}&=T_{x'w_{sr}}T_{stry_{1}}\\
&=(v_{s}-v_{s}^{-1})T_{x'w_{sr}}T_{try_{1}}+T_{x'w_{sr}s}T_{try_{1}}\\
&=(v_{s}-v_{s}^{-1})T_{x'w_{sr}}T_{try_{1}}+(v_{r}-v_{r}^{-1})T_{x'w_{sr}sr}T_{try_{1}}+T_{x'w_{sr}sr}T_{ty_{1}}.
\end{aligned}
$$
By Lemma 5.7(2), $deg\ ((v_{s}-v_{s}^{-1})T_{x'w_{sr}}T_{try_{1}})\leq L(srsr)$.
By Lemma 5.4, $deg\ ((v_{r}-v_{r}^{-1})T_{x'w_{sr}sr}T_{try_{1}})\leq L(rsr)$. If $\mathcal{L}(ty_{1})=\{t\}$,
then by Lemma 5.2, $T_{x'w_{sr}sr}T_{ty_{1}}=T_{x'w_{sr}srty_{1}}$. If $\mathcal{L}(ty_{1})=\{s,t\}$, we assume $ty_{1}=w_{st}\cdot y_{2}$ for some $y_{2}\in W$. Then
$$
\begin{aligned}
T_{x'w_{sr}sr}T_{ty_{1}}&=T_{x'w_{sr}sr}T_{stsy_{2}}\\
&=(v_{s}-v_{s}^{-1})T_{x'w_{sr}sr}T_{tsy_{2}}+T_{x'w_{sr}srs}T_{tsy_{2}}.
\end{aligned}
$$
If $\mathcal{L}(tsy_{2})=\{t\}$, then by Lemma 5.2, $T_{x'w_{sr}sr}T_{tsy_{2}}=T_{x'w_{sr}srtsy_{2}}$.
By Lemma 5.4, $deg\ (T_{x'w_{sr}srs}T_{tsy_{2}})=deg\ (T_{x'w_{sr}srsr}T_{rtsy_{2}})\leq L(sr)$.
If $\mathcal{L}(tsy_{2})=\{r,t\}$, we assume $sy_{2}=w_{sr}\cdot y_{3}$ for some $y_{3}\in W$.
By Lemma 5.2, $T_{x'w_{sr}sr}T_{tsy_{2}}=T_{x'w_{sr}srt}T_{w_{sr}y_{3}}=T_{x'w_{sr}srtw_{sr}y_{3}}$.
By Lemma 5.7(1), we get $deg\ (T_{x'w_{sr}srs}T_{tsy_{2}})=deg\ (T_{x'w_{sr}srs}T_{tw_{sr}y_{3}})=deg\ (T_{x'w_{sr}srsr\cdot rt}T_{w_{sr}y_{3}})$ $\leq L(rsr)$.
\\(2) Apply the $\mathbb{Z}[v,v^{-1}]-$algebra antiautomorphism $T_{w}\longrightarrow T_{w^{-1}}$ of $\mathcal{H}$ to (1).\hfill $\square$

\medskip

\noindent{\bf Proposition 5.9.}
We have $deg\ (T_{x}T_{y})\leq N$ for all $x,y \in W$ in case 3.

\noindent{\bf Proof.}
Now we prove Theorem 2.1 for this case. We use induction on $l(y)$. When $l(y)\leq 4$, the theorem is clear by Lemma 1.4(2). Now assume $l(y)=n\geq 5$
and the theorem is true for $y'$ if $l(y')<l(y)$. We assume $x=x_{1}\cdot w$, $y=u\cdot y_{1}$, $w,u\in W_{sr}$, $\mathcal{R}(x_{1}),\ \mathcal{L}(y_{1})\subseteq\{t\}$. Then
$$
\begin{aligned}
T_{x}T_{y}&=T_{x_{1}}(T_{w}T_{u})T_{y_{1}}\\
&=T_{x_{1}}(\sum_{p\in W_{sr}}f_{w,u,p}T_{p})T_{y_{1}}\\
&=\sum_{p\in W_{sr}}f_{w,u,p}T_{x_{1}p}T_{y_{1}}.
\end{aligned}
$$

We will prove $deg\ (f_{w,u,p}T_{x_{1}p}T_{y_{1}})\leq N$ for all $p\in W_{sr}$. If $l(p)\geq 6$ or $p=srsrs$, by Lemma 1.4(2), $deg\ f_{w,u,p}\leq L(p)\leq L(w_{sr})$. By Lemma 5.2, $T_{x_{1}p}T_{y_{1}}=T_{x_{1}py_{1}}$. So we have
$deg\ (f_{w,u,p}T_{x_{1}p}T_{y_{1}})\leq N$. If $l(p)\leq 4$ or $p=rsrsr$, $l(x_{1})\leq 2$ or $l(y_{1})\leq 2$,
then $deg\ f_{w,u,p}\leq L(rsrsr)$. On the other hand, $deg\ (T_{x_{1}p}T_{y_{1}})\leq 2L(s)$. So $deg\ (f_{w,u,p}T_{x'p}T_{y'})\leq L(srsrsrs)\leq N$. From now on we assume $l(p)\leq 4$ or $p=rsrsr$, $x_{1}=x_{2}\cdot rst$, $y_{1}=tsr\cdot y_{2}$, $\mathcal{R}(x_{2}),\mathcal{L}(y_{2})\subseteq\{s\}$,
$\mathcal{R}(x_{2}r)=\mathcal{L}(ry_{2})=\{r\}$, $\mathcal{R}(x_{2}rs)=\mathcal{L}(sry_{2})=\{s\}$. We have 10 cases.
\\(1) $p=e$.

We have
$$
\begin{aligned}
f_{w,u,p}T_{x_{1}p}T_{y_{1}}&=\delta_{w,u^{-1}}T_{x_{2}rst}T_{tsry_{2}}\\
&=\delta_{w,u^{-1}}(v_{s}-v_{s}^{-1})T_{x_{2}rst}T_{sry_{2}}+\delta_{w,u^{-1}}T_{x_{2}rs}T_{sry_{2}}.
\end{aligned}
$$
By Lemma 5.3, we get $deg\ ((v_{s}-v_{s}^{-1})T_{x_{2}rst}T_{sry_{2}})\leq L(rs)$. Since $l(sry_{2})<l(y)$, by induction hypothesis, we have $deg\ (T_{x_{2}rs}T_{sry_{2}})\leq N$.
\\(2) $p=r$.

We have
$$
f_{w,u,p}T_{x_{1}p}T_{y_{1}}=(v_{s}-v_{s}^{-1})f_{w,u,r}T_{x_{2}rstr}T_{sry_{2}}+f_{w,u,r}T_{x_{2}rsr}T_{sry_{2}}.
$$
By Lemma 5.4, we get $deg\ ((v_{s}-v_{s}^{-1})f_{w,u,r}T_{x_{2}rstr}T_{sry_{2}})\leq L(srsr)$.
If $u=e$, we have $f_{w,u,r}\in \mathbb{Z}$. Since $l(sry_{2})<l(y)$, by induction hypothesis, we have $deg\ (T_{x_{2}rsr}T_{sry_{2}})\leq N$. So $deg\ (f_{w,u,r}T_{x_{2}rsr}T_{sry_{2}})\leq N$.
If $u\neq e$, we have $deg\ f_{w,u,r}\leq L(r)$. Since $(v_{r}-v_{r}^{-1})T_{x_{2}rsr}T_{sry_{2}}=T_{x_{2}rsr}T_{rsry_{2}}-T_{x_{2}rs}T_{sry_{2}}$,
$l(rsry_{2})<l(y)$, $l(sry_{2})<l(y)$, by induction hypothesis, we get $deg\ (T_{x_{2}rsr}T_{sry_{2}})\leq N-L(r)$. So $deg\ (f_{w,u,r}T_{x_{2}rsr}T_{sry_{2}})$ $\leq N$.
\\(3) $p=s$.

Since $deg\ f_{w,u,s}\leq L(s)$, we only need to prove $deg\ (T_{x_{1}s}T_{y_{1}})\leq N-L(s)$. We have
$$
\begin{aligned}
T_{x_{1}s}T_{y_{1}}&=T_{x_{2}rtst}T_{tsry_{2}}\\
&=(v_{s}-v_{s}^{-1})T_{x_{2}rtst}T_{sry_{2}}+T_{x_{2}rts}T_{sry_{2}}\\
&=(v_{s}-v_{s}^{-1})(v_{s}-v_{s}^{-1})T_{x_{2}rtst}T_{ry_{2}}+(v_{s}-v_{s}^{-1})T_{x_{2}rst}T_{ry_{2}}\\
&\ \ \ \ +(v_{s}-v_{s}^{-1})T_{x_{2}rt}T_{sry_{2}}+T_{x_{2}rt}T_{ry_{2}}\\
&=(v_{s}-v_{s}^{-1})(v_{s}-v_{s}^{-1})T_{x_{2}rtst}T_{ry_{2}}+(v_{s}-v_{s}^{-1})T_{x_{2}rst}T_{ry_{2}}\\
&\ \ \ \ +(v_{s}-v_{s}^{-1})T_{x_{2}rt}T_{sry_{2}}+(v_{r}-v_{r}^{-1})T_{x_{2}rt}T_{y_{2}}+T_{x_{2}t}T_{y_{2}}.
\end{aligned}
$$
By Lemma 5.3, $deg\ ((v_{s}-v_{s}^{-1})(v_{s}-v_{s}^{-1})T_{x_{2}rtst}T_{ry_{2}})\leq L(srs)$. By Lemma 5.4, $deg\ ((v_{s}-v_{s}^{-1})T_{x_{2}rst}T_{ry_{2}})\leq L(srs)$,
$deg\ ((v_{s}-v_{s}^{-1})T_{x_{2}rt}T_{sry_{2}})\leq L(srs)$, $deg\ ((v_{r}-v_{r}^{-1})T_{x_{2}rt}T_{y_{2}})\leq L(rsr)$.

Now we prove $deg\ (T_{x_{2}t}T_{y_{2}})\leq N-L(s)$. We assume $x_{2}t=x'\cdot a$, $y_{2}=b\cdot y'$, $a,b\in W_{st}$, $\mathcal{R}(x')$, $\mathcal{L}(y')\subseteq\{r\}$. Since $\mathcal{R}(x_{2})$, $\mathcal{L}(y_{2})\subseteq\{s\}$, we have 5 cases.
\\\textcircled{\small{1}} $a=t$ or $b=e$.

Then $x_{2}=e$ or $y_{2}=e$, so $T_{x_{2}t}T_{y_{2}}=T_{x_{2}ty_{2}}$.
\\\textcircled{\small{2}} $a=w_{st}$, $b=st$.

We have
$$
\begin{aligned}
T_{x_{2}t}T_{y_{2}}&=T_{x'w_{st}}T_{sty'}\\
&=(v_{s}-v_{s}^{-1})T_{x'w_{st}}T_{ty'}+T_{x'st}T_{ty'}\\
&=(v_{s}-v_{s}^{-1})(v_{s}-v_{s}^{-1})T_{x'w_{st}}T_{y'}+(v_{s}-v_{s}^{-1})T_{x'ts}T_{y'}\\
&\ \ \ \ +(v_{s}-v_{s}^{-1})T_{x'st}T_{y'}+T_{x's}T_{y'}.
\end{aligned}
$$
By Lemma 5.3, $deg\ ((v_{s}-v_{s}^{-1})(v_{s}-v_{s}^{-1})T_{x'w_{st}}T_{y'})\leq L(srs)$. Since $\mathcal{R}(x't)=\mathcal{L}(ty')=\{r,t\}$, by Lemma 5.4 and Lemma 5.7, we have $deg\ ((v_{s}-v_{s}^{-1})T_{x'ts}T_{y'})\leq L(srsr)$, $deg\ ((v_{s}-v_{s}^{-1})T_{x'st}T_{y'})\leq L(srsr)$. Since $(v_{s}-v_{s}^{-1})T_{x's}T_{y'}=T_{x's}T_{sy'}-T_{x'}T_{y'}$, $l(sy')<l(y)$, $l(y')<l(y)$, by induction hypothesis, we have $deg\ (T_{x's}T_{y'})\leq N-L(s)$.
\\\textcircled{\small{3}} $a=w_{st}$, $b=s$.

We have
$$
T_{x_{2}t}T_{y_{2}}=T_{x'w_{st}}T_{sy'}=(v_{s}-v_{s}^{-1})T_{x'w_{st}}T_{y'}+T_{x'st}T_{y'}.
$$
By Lemma 5.3, $deg\ ((v_{s}-v_{s}^{-1})T_{x'w_{st}}T_{y'})\leq L(sr)$. By Lemma 5.4, Lemma 5.7(2), $deg\ (T_{x'st}T_{y'})\leq L(rsr)$.
\\\textcircled{\small{4}} $a=st$, $b=st$.

We have
$$
T_{x_{2}t}T_{y_{2}}=T_{x'w_{st}}T_{ty'}=(v_{s}-v_{s}^{-1})T_{x'w_{st}}T_{y'}+T_{x'ts}T_{y'}.
$$
By Lemma 5.3, $deg\ ((v_{s}-v_{s}^{-1})T_{x'w_{st}}T_{y'})\leq L(sr)$. By Lemma 5.4, Lemma 5.7(1), $deg\ (T_{x'ts}T_{y'})\leq L(rsr)$.
\\\textcircled{\small{5}} $a=st$, $b=s$.

We have
$$
T_{x_{2}t}T_{y_{2}}=T_{x'st}T_{sy'}=T_{x'sts}T_{y'}.
$$
By Lemma 5.3, $deg\ (T_{x'sts}T_{y'})\leq L(r)$.
\\(4) $p=sr$.

We have
$$
f_{w,u,p}T_{x_{1}p}T_{y_{1}}=f_{w,u,sr}T_{x_{2}rsts}T_{rtsry_{2}}.
$$
By Lemma 5.8(2), we get $deg\ (T_{x_{2}rsts}T_{rtsry_{2}})\leq L(srsr)$.
Since $deg\ f_{w,u,sr}$ $\leq L(sr)$, we get $deg\ (f_{w,u,sr}T_{x_{2}rsts}T_{rtsry_{2}})\leq L(srsrsr)$.
\\(5) $p=rs$.

Similar to (4).
\\(6) $p=srs$.

We have
$$
\begin{aligned}
T_{x_{1}p}T_{y_{1}}&=T_{x_{2}rstsrs}T_{tsry_{2}}\\
&=T_{x_{2}rtstr}T_{stsry_{2}}\\
&=(v_{s}-v_{s}^{-1})T_{x_{2}rtsr}T_{stsry_{2}}+T_{x_{2}rtsr}T_{stry_{2}}\\
&=(v_{s}-v_{s}^{-1})T_{x_{2}rt}T_{srstsry_{2}}+T_{x_{2}rt}T_{srstry_{2}}.
\end{aligned}
$$
By Lemma 5.1(2), we know $\mathcal{L}(srstsry_{2})=\mathcal{L}(srstry_{2})=\{s\}$. Using Lemma 5.4, we get
$deg\ ((v_{s}-v_{s}^{-1})T_{x_{2}rt}T_{srstsry_{2}})\leq L(srs)$ and $deg\ (T_{x_{2}rt}T_{srstry_{2}})\leq L(sr)$.
Since $deg\ f_{w,u,srs}\leq L(srs)$, we have $deg\ (f_{w,u,srs}T_{x_{1}p}T_{y_{1}})\leq 2L(srs)$.
\\(7) $p=rsr$.

We have
$$
f_{w,u,p}T_{x_{1}p}T_{y_{1}}=f_{w,u,rsr}T_{x_{2}rstrsr}T_{tsry_{2}}=f_{w,u,rsr}T_{x_{2}rsrts}T_{trsry_{2}}.
$$
First we have $deg\ f_{w,u,rsr}\leq L(rsr)$. By Lemma 5.1(4), we know $r\notin \mathcal{R}(x_{2}rsrts)$. So we have 2 cases.
\\\textcircled{\small{1}} $\mathcal{R}(x_{2}rsrts)=\{s\}$.

By Lemma 5.4, we have $deg\ (T_{x_{2}rsrts}T_{trsry_{2}})\leq L(sr)$ and therefore $deg\ (f_{w,u,rsr}T_{x_{2}rsrts}$ $T_{trsry_{2}})\leq L(rsrsr)$.
\\\textcircled{\small{2}} $\mathcal{R}(x_{2}rsrts)=\{s,t\}$.

Then $\mathcal{R}(x_{2}rsr)=\{s,r\}$. We assume $x_{2}rsr=x_{4}\cdot w_{sr}$ for some $x_{4}\in W$, so
$$
\begin{aligned}
T_{x_{2}rsrts}T_{trsry_{2}}&=T_{x_{4}\cdot w_{sr}\cdot ts}T_{trsry_{2}}\\
&=T_{x_{4}\cdot w_{sr}s\cdot tst}T_{trsry_{2}}\\
&=(v_{s}-v_{s}^{-1})T_{x_{4}\cdot w_{sr}s\cdot tst}T_{rsry_{2}}+T_{x_{4}\cdot w_{sr}s\cdot ts}T_{rsry_{2}}.
\end{aligned}
$$
Since $\mathcal{R}(x_{4}\cdot w_{sr}s)=\{r\}$, $\mathcal{L}(rsry_{2})=\{r\}$ or $\{s,r\}$, by Lemma 5.3 and Lemma 5.5, we have $deg\ ((v_{s}-v_{s}^{-1})T_{x_{4}\cdot w_{sr}s\cdot tst}T_{rsry_{2}})\leq L(srs)$. First we assume $\mathcal{L}(rsry_{2})=\{r\}$, then
$T_{x_{4}\cdot w_{sr}s\cdot ts}T_{rsry_{2}}=T_{x_{4}\cdot w_{sr}sr\cdot rt}T_{srsry_{2}}$.
If $\mathcal{L}(srsry_{2})=\{s\}$, we have $deg\ (T_{x_{4}\cdot w_{sr}sr\cdot rt}T_{srsry_{2}})\leq L(sr)$ by Lemma 5.4. If $\mathcal{L}(srsry_{2})=\{s,r\}$, we have $deg\ (T_{x_{4}\cdot w_{sr}sr\cdot rt}T_{srsry_{2}})= L(r)$ by easy computation. Now we assume $\mathcal{L}(rsry_{2})=\{s,r\}$, then $rsry_{2}=w_{sr}\cdot y_{3}$ for some $y_{3}\in W$, we have
$$
\begin{aligned}
&\ \ \ \ T_{x_{4}\cdot w_{sr}s\cdot ts}T_{rsry_{2}}\\
&=T_{x_{4}\cdot w_{sr}s\cdot ts}T_{w_{sr}\cdot y_{3}}\\
&=(v_{s}-v_{s}^{-1})T_{x_{4}\cdot w_{sr}s\cdot t}T_{w_{sr}\cdot y_{3}}+T_{x_{4}\cdot w_{sr}s\cdot t}T_{sw_{sr}\cdot y_{3}}\\
&=(v_{s}-v_{s}^{-1})(v_{r}-v_{r}^{-1})T_{x_{4}\cdot w_{sr}sr\cdot t}T_{w_{sr}\cdot y_{3}}+(v_{s}-v_{s}^{-1})T_{x_{4}\cdot w_{sr}sr\cdot t}T_{rw_{sr}\cdot y_{3}}\\
&\ \ \ \ +(v_{r}-v_{r}^{-1})T_{x_{4}\cdot w_{sr}sr\cdot t}T_{sw_{sr}\cdot y_{3}}+T_{x_{4}\cdot w_{sr}sr\cdot t}T_{rsw_{sr}\cdot y_{3}}\\
&=(v_{s}-v_{s}^{-1})(v_{r}-v_{r}^{-1})T_{x_{4}w_{sr}srtw_{sr}y_{3}}+(v_{s}-v_{s}^{-1})T_{x_{4}w_{sr}srtrw_{sr}y_{3}}\\
&\ \ \ \ +(v_{r}-v_{r}^{-1})T_{x_{4}w_{sr}srtsw_{sr}y_{3}}+T_{x_{4}w_{sr}srtrsw_{sr}y_{3}}.
\end{aligned}
$$
Thus, if $\mathcal{R}(x_{2}rsrts)=\{s,t\}$, we have $deg\ (T_{x_{2}rsrts}T_{trsry_{2}})\leq L(srs)$. Therefore, $deg\ (f_{w,u,rsr}T_{x_{2}rsrts}T_{trsry_{2}})\leq L(srsrsr)$.
\\(8) $p=srsr$.

We have
$$
f_{w,u,p}T_{x_{1}p}T_{y_{1}}=f_{w,u,srsr}T_{x_{2}rstsrs}T_{trsry_{2}}.
$$
By Lemma 5.1(2), we know $\mathcal{R}(x_{2}rstsrs)=\mathcal{L}(sry_{2})=\{s\}$.
By Lemma 5.4, $deg\ (T_{x_{2}rstsrs}T_{trsry_{2}})\leq L(sr)$. Since $deg\ f_{w,u,srsr}\leq L(srsr)$, we have $deg\ (f_{w,u,srsr}$ $T_{x_{2}rstsrs}T_{trsry_{2}})\leq L(srsrsr)$.
\\(9) $p=rsrs$.

Similar to (8).
\\(10) $p=rsrsr$.

We have
$$
f_{w,u,p}T_{x_{1}p}T_{y_{1}}=f_{w,u,rsrsr}T_{x_{1}rsrsr}T_{y_{1}}.
$$
By Lemma 5.6, $deg\ (T_{x_{1}rsrsr}T_{y_{1}})\leq L(s)$, so $deg\ (f_{w,u,rsrsr}T_{x_{1}rsrsr}T_{y_{1}})\leq L(srsrsr)$.\hfill $\square$

\medskip

We have completed the proof of Theorem 2.1.

\medskip

\section{Some Consequences}

\medskip

In this rest of this paper, we set
$$M=\{w_{J}|J\subseteq S,\ |W_{J}|<\infty,\ L(w_{J})=N\}.$$
$$\Lambda=\{x\cdot u\cdot y|x,y\in W,\ u\in M\}.$$
Then we have the following results.

\noindent{\bf Proposition 6.1.}
If $x,y\in W$ satisfy $deg\ (T_{x}T_{y})=N$, then $x\in \Lambda$ and $y\in \Lambda$.

\medskip

\noindent{\bf Proof.}
When $W$ is a finite Coxeter group, it is clear. When $W$ is an affine Weyl group, we can get this conclusion from [Xie1, 3.1]. When $W$ has complete Coxeter graph, see [Xie2, 3.4]. The case of $m_{sr}=\infty$, $m_{st}=2$ and the case of $m_{sr}=m_{st}=\infty$ are obvious. So we only need to check the proofs from section 3 to section 5 to consider these three cases.

First of all, we consider case 1. Keeping the assumptions and notations in section 3, we have
$$
T_{x}T_{y}=\sum_{p\in W_{st}}f_{w,u,p}T_{x_{1}p}T_{y_{1}}.
$$
We should consider when the degree of $deg\ (f_{w,u,p}T_{x_{1}p}T_{y_{1}})$ can achieve $N$. By the proof in section 3,
when $l(p)\geq 2$, if $deg\ (f_{w,u,p}T_{x_{1}p}T_{y_{1}})=N$, then we have $deg\ f_{w,u,p}=N$, so $N=L(w_{st})$ and $w=u=p=w_{st}$. We get $x\in\Lambda$ and $y\in\Lambda$.
When $p=s$, we have $deg\ (f_{w,u,p}T_{x_{1}p}T_{y_{1}})<N$. When $p=e$ or $p=t$, we have $deg\ (f_{w,u,p}T_{x_{1}p}T_{y_{1}})\leq N$. By Lemma 1.4(4), we have $x\in\Lambda$ and $y\in\Lambda$ if the equality holds.

Then we consider case 2. Keeping the assumptions and notations in section 4, we have
$$
T_{x}T_{y}=\sum_{p\in W_{sr}}f_{w,u,p}T_{x'p}T_{y'}.
$$
We should consider when the degree of $deg\ (f_{w,u,p}T_{x'p}T_{y'})$ can achieve $N$. By the proof in section 4,
when $l(p)\geq 4$, if $deg\ (f_{w,u,p}T_{x'p}T_{y'})=N$, then we have $deg\ f_{w,u,p}=N$, so $N=L(w_{sr})$ and $w=u=p=w_{sr}$. We get $x\in\Lambda$ and $y\in\Lambda$.
When $l(p)=0$ or $2\leq l(p)\leq 3$, we have $deg\ (f_{w,u,p}T_{x'p}T_{y'})<N$. When $l(p)=1$, we have $deg\ (f_{w,u,p}T_{x'p}T_{y'})\leq N$. By Lemma 1.4(4), we have $x\in\Lambda$ and $y\in\Lambda$ if the equality holds.

At last, we consider case 3. Keeping the assumptions and notations in section 5, we have
$$
T_{x}T_{y}=\sum_{p\in W_{sr}}f_{w,u,p}T_{x_{1}p}T_{y_{1}}.
$$
We should consider when the degree of $deg\ (f_{w,u,p}T_{x_{1}p}T_{y_{1}})$ can achieve $N$. By the proof in section 5,
When $l(p)\geq 6$, if $deg\ (f_{w,u,p}T_{x_{1}p}T_{y_{1}})=N$, then we have $deg\ f_{w,u,p}=N$, so $w=u=p=w_{sr}$. We get $x\in\Lambda$ and $y\in\Lambda$.
When $2\leq l(p)\leq 5$, we have $deg\ (f_{w,u,p}T_{x_{1}p}T_{y_{1}})<N$.
When $0\leq l(p)\leq 1$, we have $deg\ (f_{w,u,p}T_{x_{1}p}T_{y_{1}})\leq N$. By Lemma 1.4(4), we have $x\in\Lambda$ and $y\in\Lambda$ if the equality holds.\hfill $\square$

\medskip

\noindent{\bf Proposition 6.2.}
For any $w_{J}\in M$, $q\leq w_{J}$, $x,y\in W$, $\mathcal{R}(x),\mathcal{L}(y)\subseteq S\setminus J$, we have $deg\ (T_{xq}T_{y})\leq N-L(q)$. In particular, $T_{xw_{J}}T_{y}=T_{xw_{J}y}$.

\medskip

\noindent{\bf Proof.}
We have
$$
\begin{aligned}
T_{xw_{J}}T_{w_{J}y}&=T_{x}(T_{w_{J}}T_{w_{J}})T_{y}\\
&=T_{x}(\sum_{p\in W_{J}}f_{w_{J},w_{J},p}T_{p})T_{y}\\
&=\sum_{p\in W_{J}}f_{w_{J},w_{J},p}T_{xp}T_{y}.
\end{aligned}
$$
By Lemma 1.4(3), we know $deg\ f_{w_{J},w_{J},q}=L(q)$, since $deg\ (T_{xw_{J}}T_{w_{J}y})\leq N$ by Theorem 2.1, we get $deg\ (T_{xq}T_{y})\leq N-L(q)$.\hfill $\square$

\medskip

\noindent{\bf Proposition 6.3.}
Assume $w_{J}\in M$, then
\\(1) The left cell of $W$ containing $w_{J}$ is $\{x\cdot w_{J}|x\in W\}=\{y\in W|\mathcal{R}(y)=J\}$.
\\(2) The right cell of $W$ containing $w_{J}$ is $\{w_{J}\cdot x|x\in W\}=\{y\in W|\mathcal{L}(y)=J\}$.

\medskip

\noindent{\bf Proof.}
(1) For $x\cdot w_{J}\in W$, we have $deg\ f_{xw_{J},w_{J},xw_{J}}=N$ since $deg\ f_{w_{J},w_{J},w_{J}}=N$, so $\beta_{xw_{J},w_{J},w_{J}x^{-1}}\neq 0$. By Lemma 1.8(3), we get $w_{J}\underset{L}\sim xw_{J}$.
So we have $\{x\cdot w_{J}|x\in W\}\subseteq \{x|x\underset{L}\sim w_{J}\}$. On the other hand, for any $x\in W$ with $x\underset{L}\sim w_{J}$, we have $\mathcal{R}(x)=\mathcal{R}(w_{J})=J$, so $\{x|x\underset{L}\sim w_{J}\}\subseteq \{x\cdot w_{J}|x\in W\}=\{y\in W|\mathcal{R}(y)=J\}$.
\\(2) The proof is similar to (1).\hfill $\square$

\medskip

\section{The Lowest Two-sided Cell}

\medskip

In this section, we fix an element $w_{J}\in M$ and let $c_{0}$ be the two-sided cell of $W$ containing $w_{J}$. Then we have

\medskip

\noindent{\bf Theorem 7.1.}
(1) The two-sided cell $c_{0}$ is the lowest two-sided cell of $W$.
\\(2) We have $\Lambda=\{w\in W|a(w)=N\}\subseteq c_{0}$.

\medskip

\noindent{\bf Proof.}
(1) For $w\in W$, we may assume $w=y\cdot z$, $z\in W_{J}$, $\mathcal{R}(y)\bigcap J=\varnothing$. It is clear that $y\cdot w_{J}\underset{R}{\leq}w$.
$$
\begin{aligned}
T_{w_{J}y^{-1}}T_{yw_{J}}&=T_{w_{J}}(\sum \limits_{z\in W}f_{y^{-1},y,z}T_{z})T_{w_{J}}\\
&=T_{w_{J}}T_{w_{J}}+T_{w_{J}}(\sum \limits_{z\in W\setminus\{e\}}f_{y^{-1},y,z}T_{z})T_{w_{J}}.
\end{aligned}
$$
Since $deg\ f_{w_{J},w_{J},w_{J}}=L(w_{J})$, we get $deg\ f_{w_{J}y^{-1},yw_{J},w_{J}}=L(w_{J})$.
So $\beta_{w_{J}y^{-1},yw_{J},w_{J}}\neq 0$ and $yw_{J}\underset{L}{\sim}w_{J}$. Thus $w_{J}\underset{L}{\sim}yw_{J}\underset{R}{\leq}w$.
So $w_{J}\underset{LR}{\leq}w$ for all $w\in W$. We get $c_{0}$ is the lowest two-sided cell of $W$.
\\(2) First, we prove $\Lambda=\{w\in W|a(w)=N\}$. For any $x\cdot u\cdot y\in \Lambda$, $x,y\in W$, $u\in M$,
we have $\beta_{ux^{-1},xuy,y^{-1}u}\neq 0$ since $deg\ f_{ux^{-1},xuy,uy}=N$. Using Lemma 1.8(3), we know $a(xuy)=N$, so $\Lambda\subseteq\{w\in W|a(w)=N\}$. On the other hand, if $a(w)=a(w^{-1})=N$, choose $x,y\in W$ such that $deg\ h_{x,y,w^{-1}}=N$. Then $\gamma_{x,y,w}\neq 0$. By Lemma 1.8(1)(4), $\beta_{y,w,x}=\beta_{x,y,w}\neq 0$.
So $deg\ f_{y,w,x^{-1}}=N$. Using Proposition 6.1, we get $w\in \Lambda$.

Now, we begin to prove $\{w\in W|a(w)=N\}\subseteq c_{0}$. For $w\in W$, $a(w)=a(w^{-1})=N$, there exists $x,y\in W$ such that $deg\ f_{x,y,w^{-1}}=N$. So $\beta_{x,y,w}=\gamma_{x,y,w}\neq 0$
and $w\underset{L}{\sim}x^{-1}$. We can choose $u\in W_{J}$ such that $l(yu)=l(y)+l(u)$ and $\mathcal{R}(yu)=J$, thus $yu\underset{LR}{\sim}w_{J}$.
Since $T_{x}T_{yu}=(T_{x}T_{y})T_{u}$ and $N$ is a bound for $(W,S,L)$, we have $deg\ f_{x,yu,w^{-1}u}=N$.
Thus $\beta_{x,yu,u^{-1}w}=\gamma_{x,yu,u^{-1}w}\neq 0$, $x\underset{L}{\sim}u^{-1}y^{-1}$. So $x^{-1}\underset{R}{\sim}yu$.
We get $w\underset{L}{\sim}x^{-1}\underset{R}{\sim}yu\underset{LR}{\sim}w_{J}$.
So we have $\{w\in W|a(w)=N\}\subseteq c_{0}$.\hfill $\square$

\medskip

\noindent{\bf Remark 7.2.}
We conjecture that $\Lambda=\{w\in W|a(w)=N\}=c_{0}$ for any weighted Coxeter groups of rank 3 with positive weight function. When $L=l$, it is true since all the elements in $c_{0}$ have the same a-function value. We hope to prove it in the general case.

\end{document}